\documentclass[number,secthm,seceqn,citesort,dvips]{arxbj}
\usepackage{cursive}  

\aid{0}
\volume{16}
\issue{4}
\pubyear{2010}
\firstpage{1086}
\lastpage{1113}
\doi{10.3150/10-BEJ250}

\makeatletter
\newremark{rem}[thm]{Remark}
\newcommand{\eqref}[1]{(\ref{#1})}
\newtheorem{condition}{Condition}[section]
\newtheorem{lem}{Lemma}[section]

\newcommand{\varphiup}{\varphi}
\newcommand{\psiup}{\psi}

\newtheorem{fact}{Fact}

\makeatother

\begin{document}
\begin{frontmatter}

\title{Functional CLT for sample covariance matrices}
\runtitle{Functional CLT for sample covariance matrices}

\begin{aug}
\author[a]{\fnms{Zhidong} \snm{Bai}\thanksref{a}\ead[label=e1]{baizd@nenu.edu.cn}},
\author[b]{\fnms{Xiaoying} \snm{Wang}\thanksref{b}\ead[label=e2]{wangxy022@gmail.com}}
\and
\author[c]{\fnms{Wang} \snm{Zhou}\thanksref{c}\ead[label=e3]{stazw@nus.edu.sg}\corref{}}
\runauthor{Z. Bai, X. Wang and W. Zhou}
\address[a]{KLASMOE and School of Mathematics and Statistics, Northeast
Normal University, Changchun, 130024,
P.R. China and Department of Statistics and Applied Probability,
National University of Singapore, Singapore 117546. \printead{e1}}
\address[b]{School of Mathematics and Physics North China Electric
Power University, Beijing, 102206, P.R. China. \printead{e2}}
\address[c]{Department of Statistics and Applied Probability, National
University of Singapore, Singapore 117546. \printead{e3}}
\end{aug}

\received{\smonth{1} \syear{2009}}
\revised{\smonth{8} \syear{2009}}

%
\begin{abstract}
Using Bernstein polynomial approximations, we prove the central limit
theorem for linear spectral statistics of sample covariance matrices,
indexed by a set of functions with continuous fourth order derivatives
on an open interval including $[(1-\sqrt{y})^2,(1+\sqrt{y})^2]$, the
support of the Mar\u{c}enko--Pastur law. We also derive
the explicit expressions for asymptotic mean and covariance functions.
\end{abstract}

%

\end{frontmatter}

\section{Introduction and main result}

Let $X_n=(x_{ij})_{p\times n}, 1 \leq i \leq p, 1 \leq j \leq n$, be an
observation matrix and $x_j=(x_{1j},\ldots,x_{pj})^{\mathrm{t}}$ be the $j$th
column of $X_n$. The sample covariance matrix is
then
\[
S_n=\frac{1}{n-1}\sum_{j=1}^n(x_j-\bar{x})(x_j-\bar{x})^*,
\]
where $\bar{x}=n^{-1}\sum_{j=1}^nx_j$ and $A^*$ is the complex
conjugate transpose of $A$. The sample covariance matrix plays an
important role in
multivariate analysis since it is an unbiased estimator of the
population covariance matrix and, more importantly, many statistics in
multivariate
statistical analysis (e.g., principle component analysis, factor
analysis and multivariate regression analysis) can be expressed as
functionals of the
empirical spectral distributions of sample covariance matrices. The
\textit{empirical spectral distribution} (ESD) of a symmetric (or
Hermitian, in the complex case) $p\times p$ matrix
$A$ is defined as
\[
F^A(x)=\frac{1}{p}\times\mbox{cardinal number of } \{j\dvt\lambda_j\leq
x\},
\]
where $\lambda_1,\ldots, \lambda_p$ are the eigenvalues of $A$.

Assuming that the magnitude of the dimension $p$ is proportional to the
sample size $n$, we will study a simplified version of sample
covariance matrices,
\[
B_n=\frac{1}{n}\sum_{j=1}^nx_jx_j^*=\frac{1}{n}X_nX_n^* ,
\]
since $F^{B_n}$ and $F^{S_n}$ have the same liming properties,
according to Theorem 11.43 in \cite{bai06}. We refer to \cite{bai99}
for a review of this
field.

The first success in finding the limiting spectral distribution (LSD)
of sample covariance matrices is due to to Mar\v{c}enko and
Pastur \cite{mp1967}. Subsequent work was done in \cite{gs1977,Jo1982,s1995,w1978} and \cite{y1986}, where it was
proven that
under suitable moment conditions on $x_{ij}$, with probability 1, the
ESD $F^{B_n}$ converges to the Mar\v{c}enko--Pastur (MP)
law $F_y$
with density function
\[
F_y'(x)=\frac{1}{2\curpi xy}\sqrt{(x-a)(b-x)} ,\qquad x\in[a,b] ,
\]
with point mass $1-1/y$ at the origin if $y>1$, where $a=(1-\sqrt{y})^2$ and
$b=(1+\sqrt{y})^2$; the constant $y$ is the dimension-to-sample-size ratio
index. The commonly used method to study the convergence of $F^{B_n}$
is the \textit{Stieltjes transform}, which is defined for any distribution
function $F$ by
\[
s_F(z)\triangleq\int\frac{1}{x-z}\,\mathrm{d}F(x),\qquad \Im z\neq0 .
\]

It is easy to see that $s_F(\bar{z})=\overline{s_F(z)}$, where $\bar
{z}$ denotes the conjugate of the complex number~$z$. As is known, the Stieltjes
transform of the MP law $s(z)\triangleq s_{F_y}$ is the unique solution
to the equation
%
\begin{equation}
s=\frac{1}{1-y-z-yzs}
\end{equation}
for each $z\in\mathbb{C}^+\triangleq\{z\in\mathbb{C}\dvt\Im z>0\}$ in the
set $\{s\in\mathbb{C}\dvt-(1-y)z^{-1}+ys \in\mathbb{C}^+\}$. Explicitly,
%
\begin{equation}\label{s}
s(z)=-\frac{1}{2} \biggl(\frac{1}{y}
-\frac{1}{yz}\sqrt{z^2-(1+y)z+(1-y)^2}-\frac{1-y}{yz} \biggr) .
\end{equation}
Here, and in the sequel, $\sqrt{z}$ denotes the square root of the
complex number $z$ with positive imaginary part.

Using a Berry--Esseen-type inequality established in terms of Stieltjes
transforms, Bai \cite{bai93} was able to show that the convergence rate of
$\mathbb{E}F^{B_n}$ to $F_{y_n}$ is $\mathrm{O}(n^{-5/48})$ or $\mathrm{O}(n^{-1/4})$,
according to whether $y_n$ is close to $1$ or not. In \cite{bmt1997},
Bai, Miao and
Tsay improved these rates in the case of the convergence in
probability. Later, Bai, Miao and Yao \cite{bmy2003} proved that
$F^{B_n}$ converges to
$F_{y_n}$ at a rate of $\mathrm{O}(n^{-2/5})$ in probability and $\mathrm{O}(n^{-2/5+\eta
})$ a.s. when $y_n=p/n$ is away from 1; when $y_n=p/n$ is close to 1,
both rates
are $\mathrm{O}(n^{-1/8})$. The exact convergence rate still remains unknown for
the ESD of sample covariance matrices.

Instead of studying the convergence rate directly, Bai and Silverstein
\cite{bs2004} considered the limiting distribution of the linear spectral
statistics (LSS) of the general form of sample covariance matrices,
indexed by a set of functions analytic on an open region covering the
support of
the LSD.
More precisely, let $\mathcal{D}$ denote any region including $[a,b]$
and $\mathcal{A}(\mathcal{D})$ be the set of analytic functions
on $\mathcal{D}$. Write
$G_n(x)=p[F^{B_n}(x)-F_{y_n}(x)]$.
Bai and Silverstein proved the central limit theorem (CLT) for the LSS,
\[
G_n(f)\triangleq\int^{\infty}_{-\infty}f(x)\,\mathrm{d}G_n(x) ,\qquad
f\in\mathcal{A}(\mathcal{D}) .
\]
Their result is very useful for testing large-dimensional hypotheses.
However, the analytic assumption on $f$ seems inflexible in
practical applications because in many cases of application, the kernel
functions $f$ can only be defined on the real line, instead of on the
complex plane.
On the other hand, it is proved in \cite{bai06} that the CLT of LSS
does not hold for indicator functions. Therefore, it is natural to ask
what the
weakest continuity condition is that should be imposed on the kernel
functions so that the CLT of the LSS holds. For the CLT for other types
of matrices, one
can refer to \cite{go06}.

In this paper, we consider the CLT for
\[
G_n(f)\triangleq\int^{\infty}_{-\infty}f(x)\,\mathrm{d}G_n(x) ,\qquad f\in C^4(\mathcal{U}) ,
\]
where $\mathcal{U}$ denotes any open
interval including $[a,b]$ and $C^4(\mathcal{U})$ denotes the set of functions
$f\dvtx \mathcal{U} \rightarrow\mathbb{C}$ which
have continuous fourth order derivatives.


Denote by $\underline{s}(z)$ the Stieltjes transform of $\underline
{F}_y(x)=(1-y)\mathbb{I}_{(0, \infty)}(x)+yF_y(x)$ and set
$k(z)=\underline{s}(z)/(\underline{s}(z)+1)$, where, for $x\in\mathbb
{R}$, $\underline{s}(x)=\lim_{z\rightarrow x+\mathrm{i}0}\underline{s}(z)$.

Our main result is as follows.

\begin{thm} \label{main}
Assume that:
\begin{longlist}[(a)]
\item[(a)] for each $n$, $X_n=(x_{ij})_{p\times n}$, where $x_{ij}$ are
independent identically distributed (i.i.d.) for all $i,j$ with $\mathbb
{E}x_{11}=0$,
$\mathbb{E}|x_{11}|^2=1$, $\mathbb{E}|x_{11}|^8 < \infty$ and if
$x_{ij}$ are complex variables, $\mathbb{E}x_{11}^2=0$;

\item[(b)] $y_n=p/n\rightarrow y \in(0,\infty)$ and $y\neq1$.
\end{longlist}

The LSS $G_n=\{G_n(f)\dvt f\in C^4(\mathcal{U})\}$ then converges weakly
in finite dimensions to a Gaussian process $G=\{G(f)\dvt f\in C^4(\mathcal
{U})\}$
with mean function
%
\begin{equation} \label{mean}
\mathbb{E}G(f)=\frac{\kappa_1}{2\curpi}\int^b_a f'(x)\operatorname{arg} \bigl(1-yk^2(x) \bigr)\,\mathrm{d}x-
\frac{\kappa_2}{\curpi}\int^b_af(x){\Im} \biggl(\frac{yk^3(x)}{1-yk^2(x)} \biggr)\,\mathrm{d}x
\end{equation}
and covariance function
%
\begin{eqnarray}
c(f,g)&\triangleq&
\mathbb{E}[\{G(f)-\mathbb{E}G(f)\}\{G(g)-\mathbb{E}G(g)\}]\nonumber
\\
\label{covariance1}
&=&\frac{\kappa_1+1}{2\curpi^2}\int^b_a\int^b_af'(x_1)g'(x_2)
\ln\biggl|\frac{\overline{\underline{s}(x_1)}-\underline{s}(x_2)}{\underline
{s}(x_1)-\underline{s}(x_2)} \biggr|\,\mathrm{d}x_1\,\mathrm{d}x_2
\\
\label{covariance2}
&&{}-\frac{\kappa_2y}{2\curpi^2}\int^b_a\int^b_af'(x_1)g'(x_2){\Re
}[k(x_1)k(x_2)-\overline{k(x_1)}k(x_2)]\,\mathrm{d}x_1\,\mathrm{d}x_2 ,
\end{eqnarray}
where the parameter $\kappa_1=|\mathbb{E}x_{11}^2|^2$ takes the value
$1$ if $x_{ij}$ are real, $0$ otherwise, and
$\kappa_2=\mathbb{E}|x_{11}|^4-\kappa_1-2$.
\end{thm}

\begin{rem}
In the definition of $G_n(f)$, $\theta=\int f(x)\,\mathrm{d}F(x)$ can be
regarded as a population parameter. The linear spectral statistic
$\hat{\theta}=\int f(x)\,\mathrm{d}F_n(x)$ is then an estimator of $\theta$. We
remind the reader that the center $\theta=\int f(x)\,\mathrm{d}F(x)$, rather than
$E \int f(x)\,\mathrm{d}F_n(x)$, has its strong statistical meaning
in the application of Theorem \ref{main}. Using the limiting distribution
of $G_n(f)=n(\hat{\theta}-\theta)$, one may perform a statistical test
of the ideal hypothesis. However, in this test, one cannot apply the
limiting distribution of
$n(\hat{\theta}-\mathbb{E}\hat{\theta})$, which was studied in \cite{pl}.
\end{rem}

The strategy of the proof is to use Bernstein polynomials to
approximate functions in $C^4(\mathcal{U})$. This will be done in
Section \ref{bpa}. The problem is then reduced to the analytic case.
The truncation and renormalization steps are in Section \ref{smp}.
The convergence of the empirical processes is proved in
Section \ref{con}. We derive the mean function of the limiting
process in Section \ref{mf}.

\section{Bernstein polynomial approximations}\label{bpa}

It is well known that if $\tilde f(y)$ is a continuous function on the
interval $[0,1]$,
then the Bernstein polynomials
\[
\tilde f_m(y)=\sum_{k=0}^{m}\pmatrix{m\cr k}
y^k(1-y)^{m-k}\tilde f \biggl(\frac{k}{m} \biggr)
\]
converge to $\tilde f(y)$ uniformly on $[0,1]$ as $m\rightarrow\infty$.

Suppose that $\tilde f(y)\in C^4[0,1]$. A Taylor expansion gives
\begin{eqnarray*}
\tilde f \biggl(\frac{k}{m} \biggr)
&=&\tilde f(y)+ \biggl(\frac{k}{m}-y \biggr)\tilde f'(y)+\frac{1}{2} \biggl(\frac
{k}{m}-y \biggr)^2\tilde f''(y)\\
&&{}+\frac{1}{3!} \biggl(\frac{k}{m}-y \biggr)^3\tilde f^{(3)}(y)+\frac{1}{4!} \biggl(\frac
{k}{m}-y \biggr)^4\tilde f^{(4)}(\xi_y) ,
\end{eqnarray*}
where $\xi_y$ is a number between
$k/m$ and $y$. Hence,
%
\begin{equation} \label{appg}
\tilde f_m(y)-\tilde f(y)=\frac{y(1-y)\tilde f''(y)}{2m}+\mathrm{O} \biggl(\frac
{1}{m^2} \biggr) .
\end{equation}

For the function $f\in C^4(\mathcal{U})$, there exist $0<a_{\mathrm{l}}<a< b<b_{\mathrm{r}}$
such that $[a_{\mathrm{l}},b_{\mathrm{r}}]\subset
\mathcal{U}.$
If we let $\epsilon\in(0,1/2)$ and perform a linear transformation
$y=Lx+c$, where $L=(1-2\epsilon)/(b_{\mathrm{r}}-a_{\mathrm{l}})$ and $c=((a_{\mathrm{l}}+b_{\mathrm{r}})\epsilon
-a_{\mathrm{l}})/(b_{\mathrm{r}}-a_{\mathrm{l}})$, then
$y\in[\epsilon, 1-\epsilon]$ if $x\in[a_{\mathrm{l}},b_{\mathrm{r}}]$. Define $\tilde
f(y)\triangleq f((y-c)/L)=f(x), y \in[\epsilon, 1-\epsilon]$ and
\[
f_m(x)\triangleq
\tilde f_m(y)=\sum_{k=0}^{m}\pmatrix{m\cr k}y^k(1-y)^{m-k}\tilde f \biggl(\frac
{k}{m} \biggr) .
\]

From \eqref{appg}, we have
\begin{eqnarray*}
f_m(x)-f(x)=\tilde f_m(y)-\tilde f(y)=\frac{y(1-y)\tilde f''(y)}{2m}+\mathrm{O}
\biggl(\frac{1}{m^2} \biggr) .
\end{eqnarray*}

Since $\tilde h(y)\triangleq y(1-y)\tilde f''(y)$ has a second order
derivative, we can once again use Bernstein polynomial approximation to get
\[
\tilde h_m(y)-\tilde h(y)=\sum_{k=0}^{m}
\pmatrix{m\cr k}y^k(1-y)^{m-k}\tilde h \biggl(\frac{k}{m} \biggr)-\tilde h(y)=\mathrm{O} \biggl(\frac{1}{m}
\biggr) .
\]

So, with $h_m(x)=\tilde h_m(y)$,
\[
f(x)=f_m(x)-\frac{1}{2m}h_m(x)+\mathrm{O} \biggl(\frac{1}{m^2} \biggr) .
\]

Therefore, $G_n(f)$ can be split into three parts:
\begin{eqnarray*}
G_n(f)&=&p\int_{-\infty}^{\infty}f(x)[F^{B_n}-F_{y_n}](\mathrm{d}x)\\
&=&p \int f_m(x)[F^{B_n}-F_{y_n}](\mathrm{d}x)-\frac{p}{2m}\int
h_m(x)[F^{B_n}-F_{y_n}](\mathrm{d}x) \\
&&{} +p\int \biggl(f(x)-f_m(x)+\frac{1}{2m}h_m(x) \biggr)[F^{B_n}-F_{y_n}](\mathrm{d}x)\\
&=&\Delta_1+\Delta_2+\Delta_3 .
\end{eqnarray*}

For $\Delta_3$, under the conditions in Theorem \ref{main}, by Lemma \ref{cr}
in the \hyperref[app]{Appendix},
\[
\|F^{B_n}-F_{y_n}\|=\mathrm{O}_p(n^{-2/5}) ,
\]
where $a=\mathrm{O}_p(b)$ means that $\lim_{x\to\infty}\lim_{n\to\infty}
P(|a/b| \geq x)=0$.

Taking $m^2=[n^{3/5+\epsilon_0}]$ for some $\epsilon_0>0$ and using
integration by parts, we have that
\begin{eqnarray*}
\Delta_3&=&-p\int \biggl(f(x)-f_m(x)+\frac{1}{2m}h_m(x) \biggr)' \bigl(F_n(x)-F(x) \bigr) \,\mathrm{d}x\\
&=&\mathrm{O}_p(n^{-\epsilon_0})
\end{eqnarray*}
since $ (f(x)-f_m(x)+\frac{1}{2m}h_m(x) )'=\mathrm{O}(m^{-2})$.
From now on, we choose $\epsilon_0=1/20$, so $m=[n^{13/40}].$

Note that $f_m(x)$ and $h_m(x)$ are both analytic. Based on Conditions
\ref{cond4.1} and \ref{cond4.2} in Section \ref{con} and a martingale CLT (\cite{bi1995},
Theorem 35.12), replacing $f_m$ by $h_m$, we obtain
\[
\Delta_2=\frac{\mathrm{O}(\Delta_1)}{m}=\mathrm{o}_p(1) .
\]

It suffices to consider $\Delta_1=G_n(f_m)$. Clearly, the two
polynomials $f_m(x)$ and $\tilde f_m(y)$, defined only on the real
line, can be extended to
$[a_{\mathrm{l}},b_{\mathrm{r}}]\times[-\xi,\xi]$ and $[\epsilon,1-\epsilon]\times[-L\xi
,L\xi]$, respectively.

Since $\tilde f \in C^4[0,1]$, there exists a constant $M$ such that
$|\tilde f(y)|<M\ \forall y \in[\epsilon,1-\epsilon].$ Noting that
for $(u,v)\in
[\epsilon,1-\epsilon]\times[-L\xi,L\xi]$,
\begin{eqnarray*}
|u+\mathrm{i}v|+|1-(u+\mathrm{i}v)|&=&\sqrt{u^2+v^2}+\sqrt{(1-u)^2+v^2}\\
&\leq&
u \biggl[1+\frac{v^2}{2u^2} \biggr]+(1-u) \biggl[1+\frac{v^2}{2(1-u)^2} \biggr]\leq1+\frac
{v^2}{\epsilon} ,
\end{eqnarray*}
we have, for $y=Lx+c=u+\mathrm{i}v$,
\[
|\tilde f_m(y)|= \Biggl|\sum^{m}_{k=0}\pmatrix{m\cr k}y^k(1-y)^{m-k}\tilde f
\biggl(\frac{k}{m} \biggr) \Biggr|
\leq M \biggl(1+\frac{v^2}{\epsilon} \biggr)^m.
\]

If we take $|\xi|\leq L/\sqrt{m}$, then $ |\tilde f_m(y)|\leq M
(1+L^2/(m\epsilon) )^m \rightarrow Me^{L^2/\epsilon}$ as
$m\rightarrow\infty.$ Therefore, $ \tilde f_m(y)$ is bounded when $y\in
[\epsilon,1-\epsilon]\times[-L/\sqrt{m},L/\sqrt{m}]$. In other words,
$ f_m(x)$ is
bounded when $x\in[a_{\mathrm{l}},b_{\mathrm{r}}]\times[-1/\sqrt{m},1/\sqrt{m}]$.

Let $v=1/\sqrt{m}=n^{-13/80}$ and $\gamma_m$ be the contour formed by
the boundary of the rectangle with vertices $(a_{\mathrm{l}}\pm \mathrm{i}v)$ and
$(b_{\mathrm{r}}\pm \mathrm{i}v)$. Similarly, one can show that $h_m(x)$, $f'_m(x)$ and $h'_m(x)$
are bounded on $\gamma_m$.

\section{Simplification by truncation and normalization}\label{smp}
In this section, we will truncate the variables at a suitable level and
renormalize the truncated variables. As we will see, the truncation and
renormalization do not affect the weak limit of the spectral process.

By condition (a) in Theorem \ref{main}, for any $\delta>0$,
\[
\delta^{-8}\mathbb{E}|x_{11}|^8\mathbb{I}_{\{|x_{11}|\geq\sqrt{n}\delta
\}}\rightarrow0 ,
\]
which implies the existence of a sequence $\delta_n\downarrow0$ such
that
\[
\delta_n^{-8}\mathbb{E}|x_{11}|^8\mathbb{I}_{\{
|x_{11}|\geq
\sqrt{n}\delta_n\}}\rightarrow0
\]
as $n\rightarrow\infty$. Let $\hat{x}_{ij}=x_{ij}\mathbb{I}_{\{
|x_{ij}|\leq\sqrt{n}\delta_n\}}$ and
$\tilde{x}_{ij}=(\hat{x}_{ij}-\mathbb{E}\hat{x}_{ij})/\sigma_n$, where
$\sigma_n^2=\mathbb{E}|\hat{x}_{ij}-\mathbb{E}\hat{x}_{ij}|^2$. We then have
$\mathbb{E}\tilde{x}_{ij}=0$ and $\sigma_n^2\rightarrow1$ as
$n\rightarrow\infty$. We use $\hat{X}_n$ and $\tilde{X}_n$ to denote
the analogs of
$X_n$ when the entries $x_{ij}$ are replaced by $\hat{x}_{ij}$ and
$\tilde{x}_{ij}$, respectively; let $\hat{B}_n$ and $\tilde{B}_n$ be
analogs of
$B_n$, and let $\hat{G}_n$ and $\tilde{G}_n$ be analogs of $G_n$. We
then have
\begin{eqnarray}\label{3.1}
P(G_n\neq\hat{G}_n)&\leq &P(B_n\neq\hat{B}_n)\leq np P\bigl(|x_{11}|\geq\sqrt
{n}\delta_n\bigr)\nonumber\\[-8pt]\\[-8pt]
&\leq&
pn^{-3}\delta_n^{-8}\mathbb{E}|x_{11}|^8\mathbb{I}_{\{|x_{11}|\geq\sqrt
{n}\delta_n\}}=\mathrm{o}(n^{-2}) .\nonumber
\end{eqnarray}

From Yin, Bai and Krishnaiah \cite{ybk1988}, we know that $\lambda
_{\max}^{\hat{B}_n}$ and $\lambda_{\max}^{\tilde{B}_n}$ are
a.s. bounded
by $b=(1+\sqrt{y})^2$. Let $\lambda_j^A$ denote the $j$th largest
eigenvalue of matrix $A$. Since
\[
|\sigma_n^2-1|\leq2\mathbb{E}|x_{11}|^2\mathbb{I}_{\{|x_{11}|\geq\sqrt
{n}\delta_n\}}\leq
2\bigl(\sqrt{n}\delta_n\bigr)^{-6}\mathbb{E}|x_{11}|^8\mathbb{I}_{\{|x_{11}|\geq
\sqrt{n}\delta_n\}}=\mathrm{o}(\delta_n^2n^{-3})
\]
and
\[
|\mathbb{E}\hat{x}_{11}|^2\leq\mathbb{E}|x_{11}|^2\mathbb{I}_{\{
|x_{11}|\geq\sqrt{n}\delta_n\}}\leq \mathrm{o}(\delta_n^2n^{-3}) ,
\]
we have
%
\begin{eqnarray}\label{3.2}
&& \biggl|\int f(x)\,\mathrm{d}\hat{G}_n(x)-\int f(x)\,\mathrm{d}\tilde{G}_n(x)
 \biggr|\nonumber\\
 &&\quad\leq
K\sum_{j=1}^p|\lambda_j^{\hat{B}_n}-\lambda_j^{\tilde{B}_n}|\nonumber\\
&&\quad\leq K \bigl(\operatorname{tr}(\hat{X}_n-\tilde{X}_n)(\hat{X}_n-\tilde
{X}_n)^*\bigr)^{1/2} \\
&&\quad\leq
2(1-\sigma_n^{-1})^2\operatorname{tr}\hat{B}_n+2\sigma_n^{-2}\operatorname{tr}\mathbb
{E}\hat{X}_n\mathbb{E}\hat{X}_n^*\nonumber\\ \nonumber
&&\quad\leq
\frac{2(1-\sigma_n^2)^2}{\sigma_n^2(1+\sigma_n)^2} p \lambda_{\max}^{\hat{B}_n} +2\sigma_n^{-2}np
|\mathbb{E}\hat{x}_{11}|^2=\mathrm{o}(\delta_n^2n^{-1}) .
\end{eqnarray}

From the above estimates in \eqref{3.1} and \eqref{3.2}, we obtain
\[
\int f(x)\,\mathrm{d}G_n(x)=\int f(x)\,\mathrm{d}\tilde{G}_n(x)+\mathrm{o}_p(1) .
\]

Therefore, we only need to find the limiting distribution of $\int
f(x)\,\mathrm{d}\tilde{G}_n(x)$ with the conditions that $\mathbb{E}\tilde{x}_{11}=0$,
$\mathbb{E}|\tilde{x}_{11}|^2=1$, $\mathbb{E}|\tilde{x}_{11}|^8 < \infty
$ and $\mathbb{E}\tilde{x}_{11}^2=\mathrm{o}(n^{-2})$ for complex variables. For brevity,
in the sequel, we shall suppress the superscript on the variables and
still use $x_{ij}$ to denote the truncated and renormalized variable
$\tilde{x}_{ij}$. Note that in this paper, we use $K$ as a generic
positive constant which is independent of $n$ and which may differ from
one line to the next.

\section{Convergence of $\Delta-{\mathbb{E}}\Delta$}\label{con}

If we let $\underline{B}_n=n^{-1}X_n^*X_n$, then
$F^{\underline{B}_n}(x)=(1-y_n)\mathbb{I}_{(0,\infty)}(x)+y_nF^{B_n}(x)$.
Correspondingly, we define
$\underline{F}_{y_n}(x)=(1-y_n)\mathbb{I}_{(0,\infty
)}(x)+y_nF_{y_n}(x)$. Let $s_n(z)$ and $s_n^0(z)$ be the Stieltjes
transforms of $F^{B_n}$ and
$F_{y_n}$, respectively; let $\underline{s}_n(z)$ and $\underline
{s}_n^0(z)$ be the Stieltjes transforms of $F^{\underline{B}_n}$ and
$\underline{F}_{y_n}$,
respectively. By Cauchy's theorem, we then have
\[
\Delta_1
=\frac{1}{2\curpi \mathrm{i}}\int\oint_{\gamma_m}\frac{f_m(z)}{z-x}p[F^{B_n}-F_{y_n}](\mathrm{d}x)\,\mathrm{d}z
=-\frac{1}{2\curpi\mathrm{i}}\oint_{\gamma_m}f_m(z)p[s_n(z)-s_n^0(z)]\,\mathrm{d}z .
\]

It is easy to verify that
\[
G_n(x)=p[F^{B_n}(x)-F_{y_n}(x)]=n[F^{\underline{B}_n}(x)-\underline
{F}_{y_n}(x)] .
\]

Hence, we only need to consider $y\in(0,1)$. We shall use the following
notation:
\begin{eqnarray*}
r_j&=&\bigl(1/\sqrt{n}\bigr)x_j ,\qquad D(z)=B_n-zI_p ,\qquad D_j(z)=D(z)-r_jr_j^* ,\\
\beta_j(z)&=&\frac{1}{1+r_j^*D_j^{-1}(z)r_j} ,\qquad
\bar{\beta}_j(z)=\frac{1}{1+(1/n)\operatorname{tr}D_j^{-1}(z)} ,\\
b_n(z)&=&\frac
{1}{1+(1/n)\mathbb{E}\operatorname{tr}D_j^{-1}(z)} ,\qquad
\varepsilon_j(z)=r_j^*D_j^{-1}(z)r_j-\frac{1}{n}\operatorname{tr}D_j^{-1}(z) ,\\
\delta_j(z)&=&r_j^*D_j^{-1}(z)r_j-\frac{1}{n}\operatorname{tr}\mathbb{E}D_j^{-1}(z)
\end{eqnarray*}
and equalities
%
\begin{eqnarray}
\label{4.1}
D^{-1}(z)-D^{-1}_j(z)&=&-\beta
_j(z)D_j^{-1}(z)r_jr_j^*D_j^{-1}(z) ,\\
\label{4.2}
\beta_j(z)-\bar{\beta}_j(z)&=&-\beta_j(z)\bar{\beta
}_j(z)\varepsilon_j(z)
=-\bar{\beta}^2_j(z)\varepsilon_j(z)+\beta_j(z)\bar{\beta
}^2_j(z)\varepsilon^2_j(z) ,\\
\label{4.3}
\beta_j(z)-b_n(z)&=&-\beta_j(z)b_n(z)\delta
_j(z)=-b^2_n(z)\delta_j(z)+\beta_j(z)b^2_n(z)\delta^2_j(z) .
\end{eqnarray}

Note that by (3.4) of Bai and Silverstein \cite{bs1998}, the quantities
$\beta_j(z)$, $\bar{\beta}_j(z)$ and $b_n(z)$ are bounded in absolute
value by $|z|/v$.

Denote the $\sigma$-field generated by $r_1,\ldots,r_j$ by $\mathcal
{F}_j=\sigma(r_1,\ldots,r_j)$, and let conditional expectations
$\mathbb{E}_j(\cdot)=\mathbb{E}(\cdot|\mathcal{F}_j)$ and $\mathbb
{E}_0(\cdot)=\mathbb{E}(\cdot)$. Using the equality
%
\begin{equation}\label{4.1}
D^{-1}(z)-D^{-1}_j(z)=-\beta_j(z)D_j^{-1}(z)r_jr_j^*D_j^{-1}(z) ,
\end{equation}
we have the following well-known martingale decomposition:
\begin{eqnarray*}
p[s_n(z)-\mathbb{E}s_n(z)]&=&\operatorname{tr}\bigl(D^{-1}(z)-\mathbb
{E}D^{-1}(z)\bigr)=\sum_{j=1}^n\operatorname{tr}\bigl(\mathbb{E}_{j}D^{-1}(z)-\mathbb
{E}_{j-1}D^{-1}(z)\bigr)\\
&=&\sum_{j=1}^n\operatorname{tr}(\mathbb{E}_{j}-\mathbb
{E}_{j-1})\bigl(D^{-1}(z)-D^{-1}_j(z)\bigr)\\
& =&-\sum_{j=1}^n(\mathbb{E}_{j}-\mathbb
{E}_{j-1})\beta_j(z)r_j^*D_j^{-2}(z)r_j\\
&=&-\sum_{j=1}^n(\mathbb{E}_{j}-\mathbb{E}_{j-1})\frac{\mathrm{d}\log\beta
_j(z)}{\mathrm{d}z} .
\end{eqnarray*}

Integrating by parts, we obtain
\begin{eqnarray*}
\Delta_1-{\mathbb{E}}\Delta_1
&=& \frac{1}{2\curpi
\mathrm{i}}\sum_{j=1}^n(\mathbb{E}_{j}-\mathbb{E}_{j-1}) \oint_{\gamma_m}
f'_m(z)\log\frac{\overline{\beta}_j(z)}{\beta_j(z)}\,\mathrm{d}z\\
&=& \frac{1}{2\curpi
\mathrm{i}}\sum_{j=1}^n(\mathbb{E}_{j}-\mathbb{E}_{j-1}) \oint_{\gamma_m}
f'_m(z)\log\bigl(1+\varepsilon_j(z)\overline{\beta}_j(z)\bigr)\,\mathrm{d}z .
\end{eqnarray*}

Let $R_j(z)=\log(1+\varepsilon_j(z)\overline{\beta}_j(z))-\varepsilon
_j(z)\overline{\beta}_j(z)$ and write
\begin{eqnarray}
\Delta_1-{\mathbb{E}}\Delta_1 &=& \frac{1}{2\curpi
\mathrm{i}}\sum_{j=1}^n(\mathbb{E}_{j}-\mathbb{E}_{j-1}) \oint_{\gamma_m}
f'_m(z)\bigl(\varepsilon_j(z)\overline{\beta}_j(z)+R_j(z)\bigr)\,\mathrm{d}z\nonumber
\\
\label{mh}
&=&\frac{1}{2\curpi\mathrm{i}}\sum_{j=1}^n(\mathbb{E}_{j}-\mathbb{E}_{j-1})\int_{\gamma_{mh}}
f'_m(z)[\varepsilon_j(z)\overline{\beta}_j(z)+R_j(z)]\,\mathrm{d}z
\\
\label{mv}
&&{}+\frac{1}{2\curpi
\mathrm{i}}\sum_{j=1}^n(\mathbb{E}_{j}-\mathbb{E}_{j-1})\int_{\gamma_{mv}}
f'_m(z)[\varepsilon_j(z)\overline{\beta}_j(z)+R_j(z)]\,\mathrm{d}z ,
\end{eqnarray}
where here, and in the sequel, $\gamma_{mh}$ denotes the union of the
two horizontal parts of $\gamma_m$, and $\gamma_{mv}$ the union of the
two vertical
parts.

We first prove $\mbox{\eqref{mv}}\rightarrow0$ in probability. Let $A_n=\{
a-\epsilon_1\leq\lambda^{B_n}\leq b+\epsilon_1\}$ for any $0<\epsilon
_1<a-a_{\mathrm{l}}$ and
$A_{nj}=\{a-\epsilon_1\leq\lambda^{B_{nj}}\leq b+\epsilon_1\}$, where
$B_{nj}=B_n-r_jr_ j^*$ and $\lambda^{B}$ denotes all eigenvalues of
matrix $B$. By
the interlacing theorem (see \cite{rr2001}, page 328), it follows that
$A_n\subseteq A_{nj}$. Clearly, $\mathbb{I}_{A_{nj}}$ and $r_j$ are
independent. By
Yin, Bai and Krishnaiah \cite{ybk1988} and Bai and Silverstein \cite
{bs2004}, when $y\in(0,1),$ for any $l\geq0$,
\begin{eqnarray*}
P(\lambda_{\max}^{B_n}\geq
b+\epsilon_1)&=&\mathrm{o}(n^{-l})\quad\mbox{and}\\
P(\lambda_{\min}^{B_n}\leq a-\epsilon_1)&=&\mathrm{o}(n^{-l}) .
\end{eqnarray*}
We have $P(A_n^c)=\mathrm{o}(n^{-l})$ for any $l\geq0$.

By continuity of $s(z)$, for large $n$, there exist positive constants
$M_l$ and $M_u$ such that for all $ z\in\gamma_{mv}$, $M_l\leq
|y_ns(z)|\leq
M_u$. Letting $C_{nj}=\{|\overline{\beta}_j(z)|^{-1}\mathbb
{I}_{A_{nj}}> \epsilon_2\}$, where $0<\epsilon_2<M_l/2$ and $C_n=\bigcap
_{j=1}^nC_{nj}$,
we have
\begin{eqnarray*}
P(C_n^c)&=&P\Biggl(\bigcup_{j=1}^nC_{nj}^c\Biggr)\leq\sum_{j=1}^nP(C_{nj}^c)
=\sum_{j=1}^nP\{|\overline{\beta}_j(z)|^{-1}\mathbb{I}_{A_{nj}}\leq
\epsilon_2\}\\
&\leq&\sum_{j=1}^nP \biggl\{ \biggl| \frac{1}{n}\operatorname{tr}D_j^{-1}(z)-y_ns(z) \biggr|\mathbb
{I}_{A_{nj}}\geq\epsilon_2 \biggr\}+\sum_{j=1}^nP(A_{nj}^c)\\
&\leq&\frac{1}{\epsilon_2^4}\sum_{j=1}^n\mathbb{E} \biggl|
\frac{1}{n}\operatorname{tr}D_j^{-1}(z)-y_ns(z) \biggr|^4\mathbb{I}_{A_{nj}}+nP(A_n^c)\\
&\leq&\frac{1}{\epsilon_2^4}\sum_{j=1}^n\mathrm{O}(n^{-2/5})^4+nP(A_n^c)
\leq \mathrm{O}(n^{-2/5}) ,
\end{eqnarray*}
where we have used Lemma \ref{cr}. Defining $Q_{nj}=A_{nj}\cap C_{nj}$
and $Q_n=\bigcap_{j=1}^nQ_{nj}$, it is easy to show that $Q_{nj}$ is
independent of
$r_j$ and $P(Q_n^c)\leq P(A_n^c)+P(C_n^c)\rightarrow0$ as
$n\rightarrow\infty$.  \eqref{mv} now becomes
\[
\sum_{j=1}^n(\mathbb{E}_{j}-\mathbb{E}_{j-1})\int_{\gamma_{mv}}
f'_m(z)[\varepsilon_j(z)\overline{\beta}_j(z)+R_j(z)]\mathbb
{I}_{Q_{nj}}\,\mathrm{d}z+\mathrm{o}_p(1) .
\]

From the Burkholder inequality, Lemma \ref{quadric1} and the
inequalities $|n^{-1}\operatorname{tr}D_j(z)D_j(\bar z)|\mathbb{I}_{A_{nj}}\leq
1/(a-\epsilon_1-a_{\mathrm{l}})^2$ and
$|\bar{\beta}_j(z)|\mathbb{I}_{Q_{nj}}\leq1/\epsilon_2$, we have
\begin{eqnarray*}
&&\mathbb{E} \Biggl|\sum_{j=1}^n(\mathbb{E}_{j}-\mathbb{E}_{j-1})\int_{\gamma_{mv}}
f'_m(z)[\varepsilon_j(z)\overline{\beta}_j(z)]\mathbb{I}_{Q_{nj}}\,\mathrm{d}z
\Biggr|^2\\
&&\quad\leq K \|\gamma_{mv}\|^2\sum_{j=1}^n \sup_{z\in\gamma_{mv}}
\mathbb{E}|\varepsilon_j(z)\overline{\beta}_j(z)|^2\mathbb
{I}_{Q_{nj}}\\
&&\quad\leq Kn^{-13/40}\sum_{j=1}^n \sup_{z\in\gamma_{mv}}\mathbb
{E}|\varepsilon_j(z)|^2\mathbb{I}_{A_{nj}} \leq K n^{-13/40} .
\end{eqnarray*}
By Lemma \ref{quadric1}, for $z\in\gamma_{mv}$, we have
\[
\sum_{j=1}^n P\bigl(|\varepsilon_j(z)\overline{\beta}_j(z)|\mathbb
{I}_{Q_{nj}}\geq1/2\bigr)\leq K \sum_{j=1}^n
\mathbb{E}|\varepsilon_j(z)\overline{\beta}_j(z)|^4\mathbb
{I}_{Q_{nj}}\leq K/n .
\]
From the inequality $|\log(1+x)-x|\leq Kx^2$ for $|x|<1/2$, we get
%
\begin{eqnarray} \label{j1}
&&\mathbb{E} \Biggl|\sum_{j=1}^n(\mathbb{E}_{j}-\mathbb{E}_{j-1})\int_{\gamma_{mv}}
f'_m(z)R_j(z)\mathbb{I}_{Q_{nj}\cap\{|\varepsilon_j(z)\overline{\beta
}_j(z)|<1/2\}}\,\mathrm{d}z \Biggr|^2\nonumber\\
&&\quad\leq K \|\gamma_{mv}\|^2\sum_{j=1}^n \sup_{z\in
\gamma_{mv}}\mathbb{E}|R_j(z)|^2\mathbb{I}_{Q_{nj}\cap\{
|\varepsilon_j(z)\overline{\beta}_j(z)|<1/2\}} \\
&&\quad\leq Kn^{-13/40}\sum_{j=1}^n \sup_{z\in\gamma_{mv}}\mathbb{E}|\varepsilon_j(z)|^4\mathbb{I}_{A_{nj}} \leq K n^{-53/40} .
\nonumber
\end{eqnarray}

Therefore, from the above estimates, we can conclude that  \eqref{mv}
converges to 0 in probability. Similarly, for $z\in\gamma_{mh}$, we
also have
the following estimates:
\[
\sum_{j=1}^n P\bigl(|\varepsilon_j(z)\overline{\beta}_j(z)|\geq1/2\bigr)\leq K
\sum_{j=1}^n \mathbb{E}|\varepsilon_j(z)\overline{\beta}_j(z)|^4
\]
and
%
\begin{eqnarray} \label{j2}
&&\mathbb{E} \Biggl|\sum_{j=1}^n(\mathbb{E}_{j}-\mathbb{E}_{j-1})\int_{\gamma_{mh}}
f'_m(z)R_j(z)\mathbb{I}_{\{|\varepsilon_j(z)\overline{\beta
}_j(z)|<1/2\}}\,\mathrm{d}z \Biggr|^2\nonumber\\
&&\quad\leq K \|\gamma_{mh}\|^2\sum_{j=1}^n \sup_{z\in
\gamma_{mh}}\mathbb{E}|R_j(z)|^2\mathbb{I}_{\{|\varepsilon
_j(z)\overline{\beta}_j(z)|<1/2\}}\\
&&\quad\leq K \sum_{j=1}^n
\sup_{z\in\gamma_{mh}}\mathbb{E}|\varepsilon_j(z)\overline{\beta}_j(z)|^4 .
\nonumber
\end{eqnarray}
Thus, we get
\begin{eqnarray*}
\eqref{mh}&=& -\frac{1}{2\curpi
\mathrm{i}}\sum_{j=1}^n\mathbb{E}_{j}\int_{\gamma_{mh}}f'_m(z)[\varepsilon
_j(z)\overline{\beta}_j(z)]\,\mathrm{d}z+\mathrm{o}_p(1)\\
&\triangleq&-\frac{1}{2\curpi \mathrm{i}}\sum_{j=1}^nY_{nj}+\mathrm{o}_p(1) ,
\end{eqnarray*}
where $\mathrm{o}_p(1)$ follows from \eqref{j1},  \eqref{j2} and Condition \ref{cond4.1}
below. Therefore, our goal reduces to the convergence of $\sum_{j=1}^nY_{nj}$.

Since $Y_{nj}\in\mathcal{F}_j$ and $\mathbb{E}_{j-1}Y_{nj}=0$, $\{
Y_{nj},j=1,\ldots,n\}$ is a martingale difference sequence and thus $\sum
_{j=1}^nY_{nj}$
is a sum of a martingale difference sequence. In order to apply a
martingale CLT (\cite{bi1995}, Theorem 35.12) to it, we need to check
the following
two conditions:

\begin{condition}[(Lyapunov condition)]\label{cond4.1}
\[
\sum_{j=1}^n\mathbb{E}|Y_{nj}|^4\rightarrow0 .
\]
\end{condition}

\begin{condition}[(Conditional covariance)]\label{cond4.2}
\[
-\frac{1}{4\curpi^2}\sum_{j=1}^n\mathbb{E}_{j-1}[Y_{nj}(f_m)\cdot Y_{nj}(g_m)]
\]
converges to a constant $c(f,g)$ in probability, where $f,g\in
C^4(\mathcal{U})$ and $f_m, g_m$ are their corresponding Bernstein polynomial
approximations, respectively.
\end{condition}

\begin{pf*}{Proof of Condition \ref{cond4.1}}
By Lemmas \ref{quadric3} and \ref{beta}, for any $z\in\gamma_{mh}$,
\begin{eqnarray*}
\mathbb{E}|\varepsilon_j(z)|^6
&\leq&
\frac{K}{n^6} [(\mathbb{E}|x_{11}|^4\operatorname{tr}D_j^{-1}(z)D_j^{-1*}(z))^3
+\mathbb{E}|x_{11}|^{12}\operatorname{tr}(D_j^{-1}(z)D_j^{-1*}(z))^3 ]\\
&\leq&\frac{K}{n^6v^6} [n^3+\delta_n^4n^3 ]\leq
\frac{K}{n^3v^6} .
\end{eqnarray*}
Hence, we get
\begin{eqnarray*}
\sum_{j=1}^n\mathbb{E}|Y_{nj}|^4&\leq& K \sum_{j=1}^n
\int_{\gamma_{mh}}\mathbb{E}|\varepsilon_j(z)\overline{\beta}_j(z)|^4\,\mathrm{d}z\\
&\leq& K \sum_{j=1}^n\int_{\gamma_{mh}}(\mathbb{E}|\overline{\beta
}_j(z)|^{12})^{1/3}
(\mathbb{E}|\varepsilon_j(z)|^{6})^{2/3}\,\mathrm{d}z\\
&\leq&\frac{K}{nv^4}\rightarrow0 .
\end{eqnarray*}
\upqed
\end{pf*}

\begin{pf*}{Proof of Condition \ref{cond4.2}}
Note that in Cauchy's theorem, the integral formula is independent of
the choice of contour. Hence, we have
\begin{eqnarray*}
&&-\frac{1}{4\curpi^2}\sum_{j=1}^n\mathbb{E}_{j-1}[Y_{nj}(f_m)\cdot
Y_{nj}(g_m)]\\
&&\quad=\frac{-1}{4\curpi^2}\sum_{j=1}^n\mathbb{E}_{j-1} \biggl[\int
_{\gamma_{mh}}
f'_m(z)\mathbb{E}_{j}(\varepsilon_j(z)\overline{\beta}_j(z))\,\mathrm{d}z\cdot\int
_{\gamma'_{mh}}
g'_m(z)\mathbb{E}_{j}(\varepsilon_j(z)\overline{\beta}_j(z))\,\mathrm{d}z \biggr]\\
&&\quad=\frac{-1}{4\curpi^2}\int\!\!\int_{\gamma_{mh}\times\gamma'_{mh}}f'_m(z_1)g'_m(z_2)
\sum_{j=1}^n\mathbb{E}_{j-1} [\mathbb{E}_{j}(\varepsilon
_j(z_1)\overline{\beta}_j(z_1))
\mathbb{E}_{j}(\varepsilon_j(z_2)\overline{\beta}_j(z_2)) ]\,\mathrm{d}z_1\,\mathrm{d}z_2\\
&&\quad\triangleq \frac{-1}{4\curpi^2}\int\!\!\int_{\gamma_{mh}\times\gamma
'_{mh}}f'_m(z_1)g'_m(z_2)
\Gamma_n(z_1,z_2)\,\mathrm{d}z_1\,\mathrm{d}z_2 ,
\end{eqnarray*}
where $\Gamma_n(z_1,z_2)=\sum_{j=1}^n\mathbb{E}_{j-1} [\mathbb
{E}_{j}(\varepsilon_j(z_1)\overline{\beta}_j(z_1))
\mathbb{E}_{j}(\varepsilon_j(z_2)\overline{\beta}_j(z_2)) ]
$ and $\gamma'_{m}$ is the contour formed by the rectangle with
vertices $a'_{\mathrm{l}}\pm \mathrm{i}/2\sqrt{m}$ and $b'_{\mathrm{r}}\pm \mathrm{i}/2\sqrt{m}$.
Here, $0<a_{\mathrm{l}}<a'_{\mathrm{l}}<a<b<b'_{\mathrm{r}}<b_{\mathrm{r}}$, which means that the contour $\gamma
_m$ encloses the contour $\gamma'_m$.
$\gamma'_{mh}$ is the union of the horizontal parts of $\gamma'_m$.

First, we show that
\[
\Gamma_n(z_1,z_2)-\Gamma(z_1,z_2)\stackrel{\mathrm{Pr.}}{\longrightarrow}0
\quad\mbox{uniformly on }  \gamma_{mh}\times\gamma'_{mh} ,
\]
where
\[
\Gamma(z_1,z_2)=\kappa_2yk(z_1)k(z_2)-(\kappa_1+1)\ln\frac{\underline
{s}(z_1) \underline{s}(z_2)(z_1-z_2)}{\underline{s}(z_1)-\underline
{s}(z_2)} .
\]

From Lemma \ref{beta}, for all $z \in\gamma_{mh}\cup\gamma'_{mh}$ and
any $l\geq2$,
%
\begin{eqnarray}\label{beta1}
\mathbb{E}|\overline{\beta}_j(z)-b_n(z)|^l&=&\mathbb{E}\bigl|\overline{\beta}_j(z)b_n(z)n^{-1}\bigl(\operatorname{tr}D_j(z)-\mathbb{E}\operatorname{tr}D_j(z)\bigr)\bigr|^l\nonumber\\[-8pt]\\[-8pt]
&\leq& M
\bigl(\mathbb{E}\bigl|n^{-1}\bigl(\operatorname{tr}D_j(z)-\mathbb{E}\operatorname{tr}D_j(z)\bigr)\bigr|^{2l}\bigr)^{1/2}\leq
K\bigl(\sqrt{nv}\bigr)^l .\nonumber
\end{eqnarray}

This leads to
\[
\mathbb{E} \Biggl|\Gamma_{n}(z_1,z_2)-b_n(z_1)b_n(z_2)\sum_{j=1}^n
\mathbb{E}_{j-1}(\mathbb{E}_{j}\varepsilon_j(z_1)\mathbb
{E}_{j}\varepsilon_j(z_2)) \Biggr|\leq\frac{K}{\sqrt{n}v^3}=\mathrm{O} (n^{-
{1}/{80}} ) .
\]

Thus, we need to consider
%
\begin{equation} \label{b}
b_n(z_1)b_n(z_2)\sum_{j=1}^n \mathbb{E}_{j-1}(\mathbb{E}_{j}\varepsilon
_j(z_1)\mathbb{E}_{j}\varepsilon_j(z_2)) .
\end{equation}

Let $[A]_{ii}$ denote the $(i,i)$ entry of matrix $A$. For any two
$p\times p$ non-random matrices $A$ and~$B$, we have
%
\begin{eqnarray}\label{quadric}
&&\mathbb{E}(x_1^*Ax_1-n\operatorname{tr}A)(x_1^*Bx_1-n\operatorname{tr}B) \nonumber\\
&&\quad=(\mathbb{E}|x_{11}|^4-|\mathbb{E}x_{11}^2|^2-2)\sum
_{i=1}^pa_{ii}b_{ii}+|\mathbb{E}x_{11}^2|^2\sum_{i,j}^pa_{ij}b_{ij}+\sum
_{i,j}^pa_{ij}b_{ji}\\
\nonumber
&&\quad=\kappa_2\sum_{i=1}^pa_{ii}b_{ii}+\kappa_1\operatorname{tr}AB^{\mathrm{T}}+\operatorname{tr}AB ,
\end{eqnarray}
from which \eqref{b} becomes
\begin{eqnarray*}
&&(\kappa_1+1)b_n(z_1)b_n(z_2)\frac{1}{n^2}\sum_{j=1}^n \operatorname{tr}\mathbb
{E}_jD_j^{-1}(z_1)\mathbb{E}_jD_j^{-1}(z_2)\\
&&\qquad{}+\kappa_2b_n(z_1)b_n(z_2)\frac{1}{n^2}\sum_{j=1}^n\sum_{i=1}^p\mathbb
{E}_j[D_j^{-1}(z_1)]_{ii}\mathbb{E}_j[D_j^{-1}(z_2)]_{ii}\\
&&\quad\triangleq\Gamma_{n1}(z_1,z_2)+\Gamma_{n2}(z_1,z_2) .
\end{eqnarray*}

For $\Gamma_{n2}(z_1,z_2)$, by Lemmas \ref{beta}, \ref{diagonal} and
$-zs(z)(\underline{s}(z)+1)=1$, we get
\[
\Gamma_{n2}(z_1,z_2)=\kappa_2y_nk(z_1)k(z_2)+\mathrm{o}_p(1) ,
\]
where $\mathrm{o}_p(1)$ denotes uniform convergence in probability on $\gamma
_{mh}\times\gamma'_{mh}$.

It is easy to check that $k(\bar{z})=\overline{k(z)}$ since $s(\bar
{z})=\overline{\underline{s}(z)}$. As $n\rightarrow\infty$,
$a_{\mathrm{l}}\rightarrow a$
and $b_{\mathrm{r}}\rightarrow b$, we then get
\begin{eqnarray*}
&&-\frac{1}{4\curpi^2}\int\!\!\int_{\gamma_{mh}\times\gamma'_{mh}}
f'_m(z_1)g'_m(z_2)
\Gamma_{n2}(z_1,z_2)\,\mathrm{d}z_1\,\mathrm{d}z_2\\
&&\quad=-\frac{\kappa_2y_n}{4\curpi^2}\int\!\!\int_{\gamma_{mh}\times\gamma
'_{mh}}f'_m(z_1)g'_m(z_2)
k(z_1)k(z_2)\,\mathrm{d}z_1\,\mathrm{d}z_2+\mathrm{o}_p(1)\\
&&\quad\rightarrow-\frac{\kappa_2y}{2\curpi^2}\int^b_a\int^b_af'(x_1)g'(x_2){\Re
}[k(x_1)k(x_2)-\overline{k(x_1)}k(x_2)]\,\mathrm{d}x_1\,\mathrm{d}x_2 ,
\end{eqnarray*}
which is \eqref{covariance2} in Theorem \ref{main}.

For $\Gamma_{n1}(z_1,z_2)$, we will find the limit of
%
\begin{equation} \label{t1}
b_n(z_1)b_n(z_2)\frac{1}{n^2}\sum_{j=1}^n\operatorname{tr}\mathbb
{E}_jD_j^{-1}(z_1)\mathbb{E}_jD_j^{-1}(z_2) .
\end{equation}

Let $D_{ij}(z)=D(z)-r_jr_j^*-r_ir_i^*$, $\beta
_{ij}(z)=(1+r_i^*D_{ij}^{-1}(z)r_i)^{-1}$, $b_{12}(z)=(1+\frac
{1}{n}\mathbb{E}\operatorname{tr}D_{12}^{-1}(z))^{-1}$ and
$t(z)= (z-\frac{n-1}{n}b_{12}(z) )^{-1}$. Write
\[
D_j(z)+zI_p-\frac{n-1}{n}b_{12}(z)I_p=\sum_{i\neq j}^nr_ir_i^*-\frac
{n-1}{n}b_{12}(z)I_p .
\]

Multiplying by $t(z)I_p$ on the left, $D_j^{-1}(z)$ on the right and
combining with the identity
%
\begin{equation}\label{ridj}
r_i^*D_j^{-1}(z)=\beta_{ij}(z)r_i^*D_{ij}^{-1}(z) ,
\end{equation}
we obtain
\begin{eqnarray}\label{djni}
D^{-1}_j(z)&=&-t(z)I_p+\sum_{i\neq j}^nt(z)\beta
_{ij}(z)r_ir_i^*D_{ij}^{-1}(z)-\frac{n-1}{n}b_{12}(z)t(z)D_j^{-1}(z)\nonumber\\[-8pt]\\[-8pt]
&=&-t(z)I_p+b_{12}(z)A(z)+B(z)+C(z) ,\nonumber
\end{eqnarray}
where
\[
A(z)=\sum_{i\neq j}^nt(z)(r_ir_i^*-n^{-1}I_p)D_{ij}^{-1}(z),\qquad
B(z)=\sum_{i\neq j}^nt(z)\bigl(\beta_{ij}(z)-b_{12}(z)\bigr)r_ir_i^*D_{ij}^{-1}(z)
\]
and
\[
C(z)=\frac{1}{n}t(z)b_{12}(z)\sum_{i\neq
j}^n\bigl(D_{ij}^{-1}(z)-D_j^{-1}(z)\bigr) .
\]

It is easy to verify that for all $z\in\gamma_{mh}\cup\gamma'_{mh}$,
\begin{eqnarray*}
|t(z)|&= &\biggl|z+\frac{n-1}{n}\frac{1}{1+n^{-1}\mathbb{E}\operatorname{tr}D_{12}^{-1}(z)} \biggr|^{-1}
= \biggl|\frac{1+n^{-1}\mathbb{E}\operatorname{tr}D_{12}^{-1}(z)}{z(1+n^{-1}\mathbb
{E}\operatorname{tr}D_{12}^{-1}(z))+(n-1)/n} \biggr|\\
&\leq&\frac{1}{|z|} \biggl[1+\frac{1}{\Im z(1+n^{-1}\mathbb
{E}\operatorname{tr}D_{12}^{-1}(z))} \biggr]
\leq\frac{K}{v}
\end{eqnarray*}
since $a_{\mathrm{l}}\leq|z|\leq b_{\mathrm{r}}+1$. Thus, by Lemmas \ref{beta}, \ref
{quadric2} and the Cauchy--Schwarz inequality, we have
%
\begin{eqnarray} \label{j3}
\mathbb{E}|\operatorname{tr}(B(z_1)\mathbb{E}_jD_j^{-1}(z_2))|
&=&\mathbb{E} \Biggl|\sum_{i\neq
j}^nt(z_1)\bigl(\beta_{ij}(z_1)-b_{12}(z_1)\bigr)r_i^*D_{ij}^{-1}(z_1)\mathbb
{E}_jD_j^{-1}(z_2)r_i \Biggr|\nonumber\\
&\leq&\frac{Kn}{v}\mathbb{E}\bigl|\bigl(\beta
_{ij}(z_1)-b_{12}(z_1)\bigr)r_i^*D_{ij}^{-1}(z_1)\mathbb
{E}_jD_j^{-1}(z_2)r_i\bigr| \\
&\leq&\frac{Kn}{v}\frac{1}{\sqrt{n}v}\frac{1}{\sqrt{n}v^2}=\frac
{K}{v^4} .\nonumber
\end{eqnarray}

From Lemma 2.10 of Bai and Silverstein \cite{bs1998}, for any $n\times
n$ matrix $A$,
%
\begin{equation}\label{d-dj}
\bigl|\operatorname{tr}\bigl(D^{-1}(z)-D_j^{-1}(z)\bigr)A\bigr|\leq\frac{\|A\|}{v},
\end{equation}
which, combined with Lemma \ref{beta}, gives
\begin{eqnarray} \label{j4}
&&\mathbb{E}|\operatorname{tr}(C(z_1)\mathbb{E}_jD_j^{-1}(z_2))|\nonumber\\
&&\quad=\mathbb{E} \Biggl|\frac{1}{n}t(z_1)b_{12}(z_1)\sum_{i\neq
j}^n\operatorname{tr}\bigl(\bigl(D_{ij}^{-1}(z_1)-D_j^{-1}(z_1)\bigr)\mathbb{E}_jD_j^{-1}(z_2)\bigr)
\Biggr|\nonumber\\[-8pt]\\[-8pt]
&&\quad\leq\frac{K}{v}(\mathbb{E}|b_{12}(z_1)|^2)^{1/2}\bigl(\mathbb
{E}\bigl|\operatorname{tr}\bigl(D_{ij}^{-1}(z_1)-D_j^{-1}(z_1)\bigr)\mathbb
{E}_jD_j^{-1}(z_2)\bigr|^2\bigr)^{{1}/{2}}\nonumber\\
&&\quad\leq\frac{K}{v}\frac{1}{v^2}=\frac{K}{v^3} .\nonumber
\end{eqnarray}

From the above estimates \eqref{j3} and \eqref{j4}, we arrive at
%
\begin{eqnarray}\label{trdd}
&&\operatorname{tr}\mathbb{E}_jD_j^{-1}(z_1)\mathbb{E}_jD_j^{-1}(z_2)\nonumber\\[-8pt]\\[-8pt]
&&\quad=-t(z_1)\operatorname{tr}\mathbb{E}_jD_j^{-1}(z_2)+b_{12}(z_1)\operatorname{tr}\mathbb{E}_jA(z_1)\mathbb{E}_jD_j^{-1}(z_2)+\frac{K}{v^4} .\nonumber
\end{eqnarray}

Using the identity
\[
D_j^{-1}(z_2)-D_{ij}^{-1}(z_2)=-\beta
_{ij}(z_2)D_{ij}^{-1}(z_2)r_ir_i^*D_{ij}^{-1}(z_2) ,
\]
we can write
\begin{equation}\label{trad}
\operatorname{tr}\mathbb
{E}_j(A(z_1))D_j^{-1}(z_2)=A_1(z_1,z_2)+A_2(z_1,z_2)+A_3(z_1,z_2) ,
\end{equation}
where
\begin{eqnarray*}
A_1(z_1,z_2)&=&-\operatorname{tr}\sum_{i<j}t(z_1)r_ir_i^*\mathbb
{E}_j(D_{ij}^{-1}(z_1))\bigl(D_j^{-1}(z_2)-D_{ij}^{-1}(z_2)\bigr)\\
&=&-\sum_{i<j}t(z_1)\beta_{ij}(z_2)r_i^*\mathbb
{E}_j(D_{ij}^{-1}(z_1))D_{ij}^{-1}(z_2)r_ir_i^*D_{ij}^{-1}(z_2)r_i ,\\
A_2(z_1,z_2)&=&-\operatorname{tr}\sum_{i<j}^nt(z_1)\frac{1}{n}\mathbb
{E}_j(D_{ij}^{-1}(z_1))\bigl(D_j^{-1}(z_2)-D_{ij}^{-1}(z_2)\bigr)
\end{eqnarray*}
and
\[
A_3(z_1,z_2)=-\operatorname{tr}\sum_{i<j}t(z_1) \biggl(r_ir_i^*-\frac{1}{n}I_p
\biggr)\mathbb{E}_j(D_{ij}^{-1}(z_1))D_{ij}^{-1}(z_2) .
\]
From \eqref{d-dj}, we get
\begin{eqnarray}\label{a2}
|A_2(z_1,z_2)|&=& \biggl|\frac{1}{n}\sum_{i<j}t(z_1)\operatorname{tr}\bigl(D_j^{-1}(z_2)-D_{ij}^{-1}(z_2)\bigr)\mathbb{E}_jD_{ij}^{-1}(z_1)
\biggr|\nonumber\\[-8pt]\\[-8pt]
&\leq&\frac{j-1}{n}\frac{1}{v}\frac{K}{v^2}\leq\frac{K}{v^3}\nonumber
\end{eqnarray}
and by Lemma \ref{quadric1}, we have
\begin{eqnarray}\label{a3}
\mathbb{E}|A_3(z_1,z_2)|&\leq&\frac{K(j-1)}{v}\mathbb{E} \biggl|\operatorname{tr}
\biggl(r_ir_i^*-\frac{1}{n}I_p \biggr)\mathbb
{E}_j(D_{ij}^{-1}(z_1))D_{ij}^{-1}(z_2)\biggr|\nonumber\\[-8pt]\\[-8pt]
&\leq&\frac{Kn}{v}\frac{1}{\sqrt{n}v^2}=\frac{K\sqrt{n}}{v^3} .\nonumber
\end{eqnarray}

For $A_1(z_1,z_2)$, by Lemmas \ref{quadric2} and \ref{quadric3},
\begin{eqnarray*}
&&\mathbb{E} \biggl|r_i^*\mathbb
{E}_j(D_{ij}^{-1}(z_1))D_{ij}^{-1}(z_2)r_ir_i^*D_{ij}^{-1}(z_2)r_i\\
&&\qquad{}-\frac{1}{n^2}\operatorname{tr}[\mathbb
{E}_j(D_{ij}^{-1}(z_1))D_{ij}^{-1}(z_2)]\operatorname{tr}D_{ij}^{-1}(z_2) \biggr|\\
&&\quad\leq\mathbb{E} \biggl| \biggl[r_i^*\mathbb{E}_j(D_{ij}^{-1}(z_1))D_{ij}^{-1}(z_2)r_i-
\frac{1}{n}\operatorname{tr}(\mathbb{E}_j(D_{ij}^{-1}(z_1))D_{ij}^{-1}(z_2))
\biggr]r_i^*D_{ij}^{-1}(z_2)r_i \biggr|\\
&&\qquad{}+\mathbb{E} \biggl|\frac{1}{n}\operatorname{tr}(\mathbb
{E}_j(D_{ij}^{-1}(z_1))D_{ij}^{-1}(z_2))
\biggl[r_i^*D_{ij}^{-1}(z_2)r_i-\frac{1}{n}\operatorname{tr}D_{ij}^{-1}(z_2) \biggr] \biggr|
\leq\frac{K}{\sqrt{n}v^3} .
\end{eqnarray*}

Let $\varphiup_j(z_1,z_2)=\operatorname{tr}(\mathbb
{E}_j(D_j^{-1}(z_1))D_j^{-1}(z_2))$. Using the identity \eqref{d-dj},
we have
\[
|\operatorname{tr}(\mathbb{E}_j(D_{ij}^{-1}(z_1))D_{ij}^{-1}(z_2))\operatorname{tr}D_{ij}^{-1}(z_2)
-\varphiup_j(z_1,z_2)\operatorname{tr}D_j^{-1}(z_2) |\leq Knv^{-3}.
\]

Thus, in conjunction with Lemma \ref{beta}, we can get
\begin{equation}\label{a1}
\mathbb{E} \biggl|A_1(z_1,z_2) +\frac{j-1}{n^2}t(z_1)b_{12}(z_2)\varphiup
_j(z_1,z_2)\operatorname{tr}D_j^{-1}(z_2) \biggr|\leq\frac{K}{\sqrt{n}v^3} .
\end{equation}

Therefore, from \eqref{djni}--\eqref{a1}, it follows that
\begin{eqnarray*}
&&\varphiup_j(z_1,z_2)
\biggl[1+\frac{j-1}{n^2}t(z_1)b_{12}(z_1)b_{12}(z_2)\operatorname{tr}D_j^{-1}(z_2) \biggr]\\
&&\quad=-\operatorname{tr}(t(z_1)\operatorname{tr}D_j^{-1}(z_2))+A_4(z_1,z_2) ,
\end{eqnarray*}
where
$\mathbb{E}|A_4(z_1,z_2)|\leq K\sqrt{n}/v^3$.

Using Lemma \ref{beta}, the expression for $D_j^{-1}(z_2)$ in \eqref
{djni} and the estimate
\begin{eqnarray*}
\mathbb{E}|\operatorname{tr}A(z)|&=&\mathbb{E} \Biggl|\operatorname{tr}\sum_{i\neq
j}^nt(z)(r_ir_i^*-n^{-1}I_p)D_{ij}^{-1}(z) \Biggr|\\
&\leq&\frac{Kn}{v}\mathbb{E}|r_iD_{ij}^{-1}(z)r_i^*-n^{-1}\operatorname{tr}D_{ij}^{-1}(z)|
\leq\frac{K\sqrt{n}}{v^2},
\end{eqnarray*}
we find
that
\begin{eqnarray*}
&&\varphiup_j(z_1,z_2)
\biggl[1-\frac{(j-1)p}{n^2}t(z_1)b_{12}(z_1)t(z_2)b_{12}(z_2) \biggr]\\
&&\quad=-pt(z_1)t(z_2)+A_5(z_1,z_2) ,
\end{eqnarray*}
where
\[
\mathbb{E}|A_5(z_1,z_2)|\leq\frac{K\sqrt{n}}{v^3} .
\]

By Lemma \ref{beta}, we can write
\begin{eqnarray*}
&&\varphiup_j(z_1,z_2)
\biggl[1-\frac{(j-1)p}{n^2}\frac{\underline{s}_n^0(z_1)\underline
{s}_n^0(z_2)}{(\underline{s}_n^0(z_1)+1)(\underline{s}_n^0(z_2)+1)} \biggr]\\
&&\qquad=\frac{p}{z_1z_2}\frac{1}{(\underline{s}_n^0(z_1)+1)(\underline
{s}_n^0(z_2)+1)}+A_6(z_1,z_2) ,
\end{eqnarray*}
where
$\mathbb{E}|A_6(z_1,z_2)|\leq K\sqrt{n}/v^3$.

Let
\[
a_n(z_1,z_2)=\frac{y_n\underline{s}_n^0(z_1)\underline{s}_n^0(z_2)}
{(\underline{s}_n^0(z_1)+1)(\underline{s}_n^0(z_2)+1)} .
\]
\eqref{t1} can be written as
\[
a_n(z_1,z_2)\frac{1}{n}\sum_{j=1}^n \biggl(1-\frac{j-1}{n}a_n(z_1,z_2)
\biggr)^{-1}+A_7(z_1,z_2) ,
\]
where
\[
\mathbb{E}|A_7(z_1,z_2)|\leq\frac{K}{\sqrt{n}v^3} .
\]

Since
\[
a_n(z_1,z_2)\rightarrow a(z_1,z_2)=\frac{y\underline{s}(z_1)\underline
{s}(z_2)}{(\underline{s}(z_1)+1)(\underline{s}(z_2)+1)}
\]
as $n\rightarrow\infty$, we arrive at
\begin{eqnarray*}
\eqref{t1}&\stackrel{\mathrm{Pr.}}{\longrightarrow} & a(z_1,z_2)\int_0^1\frac
{1}{1-ta(z_1,z_2)}\,\mathrm{d}t=-\ln\bigl(1-a(z_1,z_2)\bigr)=-\ln\frac{l(z_1,z_2)}{\underline{s}(z_1)-\underline{s}(z_2)} ,
\end{eqnarray*}
where $l(z_1,z_2)=\underline{s}(z_1)\underline{s}(z_2)(z_1-z_2)$, which
implies that
\begin{eqnarray*}
\Gamma_{n1}(z_1,z_2)&=&(\kappa_1+1)b_n(z_1)b_n(z_2)\frac{1}{n^2}\sum
_{j=1}^n\operatorname{tr}\mathbb{E}_jD_j^{-1}(z_1)\mathbb{E}_jD_j^{-1}(z_2)\\
&=&-(\kappa_1+1)\ln(l(z_1,z_2))+(\kappa_1+1)\ln\bigl(\underline
{s}(z_1)-\underline{s}(z_2)\bigr)+\mathrm{o}_p(1) .
\end{eqnarray*}

Thus, adding the vertical parts of both contours and using the fact
that $f'_m(z)$ and $g'_m(z)$ are analytic functions, the integral of
the first term
of $\Gamma_{n1}(z_1,z_2)$ is
\begin{eqnarray*}
&&-\frac{1}{4\curpi^2}\int\!\!\int_{\gamma_{mh}\times\gamma'_{mh}}f'_m(z_1)g'_m(z_2)
(\kappa_1+1)\ln(l(z_1,z_2))\,\mathrm{d}z_1\,\mathrm{d}z_2\\
&&\quad=-\frac{\kappa_1+1}{4\curpi^2}\oint\oint_{\gamma_{m}\times\gamma'_{m}}f'_m(z_1)g'_m(z_2)
\ln(l(z_1,z_2))\,\mathrm{d}z_1\,\mathrm{d}z_2+\mathrm{O}(v)\\
&&\quad=\mathrm{o}(1) .
\end{eqnarray*}

For the second term of $\Gamma_{n1}(z_1,z_2)$, since $s(\bar
{z})=\overline{\underline{s}(z)}$, as $n\rightarrow\infty$,
$a_{\mathrm{l}}\rightarrow a$ and
$b_{\mathrm{r}}\rightarrow b$, we get
\begin{eqnarray*}
&&-\frac{\kappa_1+1}{4\curpi^2}\oint\oint_{\gamma_{m}\times\gamma'_{m}}f'_m(z_1)g'_m(z_2)
\ln\bigl(\underline{s}(z_1)-\underline{s}(z_2)\bigr)\,\mathrm{d}z_1\,\mathrm{d}z_2+\mathrm{o}_p(1)\\
&&\quad\longrightarrow\frac{\kappa_1+1}{2\curpi^2}\int^b_a\int^b_af'(x_1)g'(x_2)
\ln\biggl|\frac{\overline{\underline{s}(x_1)}-\underline{s}(x_2)}{\underline
{s}(x_1)-\underline{s}(x_2)} \biggr|\,\mathrm{d}x_1\,\mathrm{d}x_2 ,
\end{eqnarray*}
which is \eqref{covariance1} in Theorem \ref{main}.
\end{pf*}

\section{Mean function} \label{mf}

In this section, we will find the limit of
\[
\mathbb{E}G_n(f_m)=-\frac{1}{2\curpi
\mathrm{i}}\oint_{\gamma_m}f_m(z)p[\mathbb{E}s_n(z)-s_n^0(z)]\,\mathrm{d}z .
\]

We shall first consider
$M_n(z)=p[\mathbb{E}s_n(z)-s^0_n(z)]=n[\mathbb{E}\underline
{s}_n(z)-\underline{s}_n^0(z)].$

Since $D(z)+zI=\sum_{j=1}^nr_jr_j^*$, multiplying by $D^{-1}(z)$ on the
right-hand side and using \eqref{ridj}, we find that
\[
I+zD^{-1}(z)=\sum_{j=1}^nr_jr_j^*D^{-1}(z)=\sum_{j=1}^n\frac
{r_jr_j^*D_j^{-1}(z)}{1+r_j^*D_j^{-1}(z)r_j} .
\]

Taking trace, dividing by $n$ on both sides and combining with the
identity $z\underline{s}_n(z)=-1+y_n+y_nzs_n(z)$ leads to
\begin{equation}\label{snbeta}
\underline{s}_n(z)=-\frac{1}{nz}\sum_{j=1}^n\frac
{1}{1+r_j^*D_j^{-1}(z)r_j} =-\frac{1}{nz}\sum_{j=1}^n\beta_j(z) .
\end{equation}

Then, once again using \eqref{ridj} and
$A^{-1}-B^{-1}=-A^{-1}(A-B)B^{-1}$, we get
\begin{eqnarray*}
\frac{I_p}{z(\mathbb{E}\underline{s}_n(z)+1)}-D^{-1}(z)
&=&-\frac{1}{z(\mathbb{E}\underline{s}_n(z)+1)} \Biggl[\sum
_{j=1}^nr_jr_j^*+z\mathbb{E}\underline{s}_n(z) \Biggr]D^{-1}(z)\\
&=&-\frac{1}{z(\mathbb{E}\underline{s}_n(z)+1)}\sum_{j=1}^n \biggl[\beta
_j(z)r_jr_j^*D_j^{-1}(z)
-\mathbb{E}(\beta_j(z))\frac{1}{n}D^{-1}(z) \biggr].
\end{eqnarray*}

Taking trace, dividing by $p$ and taking expectation, we find that
\begin{eqnarray}\label{omga}
\omega_n(z)&=& -\frac{1}{z(\mathbb{E}\underline{s}_n(z)+1)}-\mathbb
{E}s_n(z)\nonumber\\
&=&-\frac{1}{p
z(\mathbb{E}\underline{s}_n(z)+1)}\sum_{j=1}^n\mathbb{E}(\beta
_j(z)d_j(z))\\
&\triangleq&-\frac{1}{p
z(\mathbb{E}\underline{s}_n(z)+1)}J_n(z) ,\nonumber
\end{eqnarray}
where
\[
d_j(z)=r_j^*D_j^{-1}(z)r_j-\frac{1}{n}\operatorname{tr}\mathbb{E}D^{-1}(z) .
\]

On the other hand, by the identity $\mathbb{E}\underline
{s}_n(z)=-(1-y_n)z^{-1}+y_n\mathbb{E}s_n(z)$, we have
\[
\omega_n(z)=\frac{\mathbb{E}\underline{s}_n(z)}{y_nz} \biggl(-z
-\frac{1}{\mathbb{E}\underline{s}_n(z)}+\frac{y_n}{\mathbb{E}\underline
{s}_n(z)+1} \biggr)
\triangleq\frac{\mathbb{E}\underline{s}_n(z)}{y_nz}R_n(z) ,
\]
where
\[
R_n(z)=-z-\frac{1}{\mathbb{E}\underline{s}_n(z)}+\frac{y_n}{\mathbb
{E}\underline{s}_n(z)+1} ,
\]
which implies that
%
\begin{equation}\label{sn}
\mathbb{E}\underline{s}_n(z)= \biggl(-z+\frac{y_n}{\mathbb{E}\underline{s}_n(z)+1}-R_n(z) \biggr)^{-1}.
\end{equation}

For $s_n^0(z)$, since $s_n^0(z)=(1-y_n-y_nzs_n^0(z)-z)^{-1}$ and
$\underline{s}_n^0(z)=-(1-y_n)z^{-1}+y_ns_n^0(z)$, we have
%
\begin{equation}\label{s0}
\underline{s}_n^0(z)= \biggl(-z+\frac{y_n}{\underline{s}_n^0(z)+1} \biggr)^{-1}.
\end{equation}

By \eqref{sn} and \eqref{s0}, we get
\begin{eqnarray*}
\mathbb{E}\underline{s}_n(z)-\underline{s}_n^0(z)
&=& \biggl(-z+\frac{y_n}{\mathbb{E}\underline{s}_n(z)+1}-R_n(z) \biggr)^{-1}- \biggl(-z+\frac
{y_n}{\underline{s}_n^0(z)+1} \biggr)^{-1}\\
&=&\mathbb{E}\underline{s}_n(z)\underline{s}_n^0(z) \biggl(\frac
{y_n}{\underline{s}_n^0(z)+1}-\frac{y_n}{\mathbb{E}\underline{s}_n(z)+1}+R_n(z) \biggr)\\
&=&\frac{y_n\mathbb{E}\underline{s}_n(z)\underline
{s}_n^0(z)}{(\underline{s}_n^0(z)+1)(\mathbb{E}\underline{s}_n(z)+1)}
\bigl(\mathbb{E}\underline{s}_n(z)-\underline{s}_n^0(z) \bigr)
+\mathbb{E}\underline{s}_n(z)\underline{s}_n^0(z)R_n(z) ,
\end{eqnarray*}
which, combined with \eqref{omga}, leads to
\begin{eqnarray} \label{diff}
&&n\bigl(\mathbb{E}\underline{s}_n(z)-\underline{s}_n^0(z)\bigr) \biggl(1-
\frac{y_n\mathbb{E}\underline{s}_n(z)\underline{s}_n^0(z)}{(\underline
{s}_n^0(z)+1)(\mathbb{E}\underline{s}_n(z)+1)} \biggr)\nonumber\\
&&\quad=n\mathbb{E}\underline{s}_n(z)\underline{s}_n^0(z)R_n(z)\nonumber\\[-8pt]\\[-8pt]
\nonumber
&&\quad=n\mathbb{E}\underline{s}_n(z)\underline{s}_n^0(z)\frac
{y_nz}{\mathbb{E}\underline{s}_n(z)}\omega_n(z)\\
&&\quad=-\frac{\underline{s}_n^0(z)}{\mathbb{E}\underline{s}_n(z)+1}J_n(z) .\nonumber
\end{eqnarray}

Thus, in order to find the limit of $M_n(z)=n[\mathbb{E}\underline
{s}_n(z)-\underline{s}_n^0(z)]$, it suffices to find the limit of
$J_n(z)$. Let
$\bar{d}_j(z)= r_j^*D_j^{-1}(z)r_j-\frac{1}{n}\operatorname{tr}D^{-1}(z)$ and $\bar
{J}_n(z)=\sum_{j=1}^n\mathbb{E}(\beta_j(z)\bar{d}_j(z))$. By \eqref
{4.3}, we have
\begin{eqnarray*}
J_n(z)&=&\bar{J}_n(z)+\sum_{j=1}^n\mathbb{E} \biggl[\beta_j(z) \biggl(\frac
{1}{n}\operatorname{tr}D^{-1}(z)-\frac{1}{n}\operatorname{tr}\mathbb{E}D^{-1}(z) \biggr) \biggr]\\
&=&\bar{J}_n(z)+\sum_{j=1}^n\mathbb{E} \biggl[\bigl(\beta_j(z)-b_n(z)\bigr) \biggl(\frac
{1}{n}\operatorname{tr}D^{-1}(z)-\frac{1}{n}\operatorname{tr}\mathbb{E}D^{-1}(z) \biggr) \biggr]\\
&=&\bar{J}_n(z)-T_1(z)+T_2(z) ,
\end{eqnarray*}
where, from \eqref{d-dj},
\begin{eqnarray*}
T_1(z)&=&\sum_{j=1}^n\mathbb{E} \biggl[b^2_n(z)\delta_j(z)
\biggl(\frac{1}{n}\operatorname{tr}D^{-1}(z)-\frac{1}{n}\operatorname{tr}\mathbb{E}D^{-1}(z) \biggr) \biggr]\\
&=&\sum_{j=1}^n\mathbb{E} \biggl[b^2_n(z)\delta_j(z)
\frac{1}{n} \bigl(\operatorname{tr}\bigl(D^{-1}(z)-D_j^{-1}(z)\bigr)-\operatorname{tr}\mathbb
{E}\bigl(D^{-1}(z)-D_j^{-1}(z)\bigr) \bigr) \biggr]\\
&\leq&\sum_{j=1}^n\frac{K}{\sqrt{n}v}\cdot\frac{K}{nv}=\frac{K}{\sqrt
{n}v^2} .
\end{eqnarray*}

It follows from Bai and Silverstein \cite{bs1998}, (4.3) that for $l\geq2$,
\begin{equation}\label{ded}
\mathbb{E} \biggl|\frac{1}{n}\operatorname{tr}D^{-1}(z)-\frac{1}{n}\operatorname{tr}\mathbb{E}D^{-1}(z) \biggr|^l\leq\frac{K_l}{(\sqrt{n}v)^l} .
\end{equation}
Hence,
\[
T_2(z)=\sum_{j=1}^n\mathbb{E} \biggl[\beta_j(z)b^2_n(z)\delta_j^2(z) \biggl(\frac
{1}{n}\operatorname{tr}D^{-1}(z)-\frac{1}{n}\operatorname{tr}\mathbb{E}D^{-1}(z) \biggr) \biggr]
\leq\frac{K}{\sqrt{n}v^3} .
\]

From the above estimates on $T_1$ and $T_2$, we conclude that
\[
J_n(z)=\bar{J}_n(z)+\bar{\epsilon}_n ,
\]
where here, and in the sequel, $\bar{\epsilon}_n=\mathrm{O}((\sqrt{n}v^3)^{-1})$.

We now only need to consider the limit of $\bar{J}_n(z)$. By \eqref
{4.2}, we write
\begin{eqnarray*}
\bar{J}_n(z)
&=&\sum_{j=1}^n\mathbb{E}\bigl[\bigl(\beta_j(z)-\bar{\beta}_j(z)\bigr)\varepsilon
_j(z)\bigr]+\frac{1}{n}\sum_{j=1}^n\mathbb{E}\bigl[\beta_j(z)\operatorname{tr}\bigl(D_j^{-1}(z)-D^{-1}(z)\bigr)\bigr]\\
&=&-\sum_{j=1}^n\mathbb{E}(\bar{\beta}^2_j(z)\varepsilon^2_j(z))+\sum
_{j=1}^n\mathbb{E}(\bar{\beta}^2_j(z)\beta_j(z)\varepsilon^3_j(z))
+\frac{1}{n}\sum_{j=1}^n\mathbb{E}(\beta^2_j(z)r_j^*D^{-2}_j(z)r_j)\\
&\triangleq& \bar{J}_{n1}(z)+\bar{J}_{n2}(z)+\bar{J}_{n3}(z) .
\end{eqnarray*}

From Lemmas \ref{quadric1} and \ref{beta}, we find that
\[
|\bar{J}_{n2}(z)|\leq K\sum_{j=1}^n (\mathbb{E}|\varepsilon^6_j(z)|
)^{1/2}\leq\frac{K}{\sqrt{n}v^3} .
\]

By Lemma \ref{beta}, $\bar{\beta}_j(z),
\beta_j(z)$ and $b_n(z)$ can be replaced by $-z\underline{s}(z)$, and
so we get
\[
\bar{J}_{n3}(z)=z^2\underline{s}^2(z)\frac{1}{n^2}\sum_{j=1}^n\mathbb
{E}\operatorname{tr}D^{-2}_j(z)+\bar{\epsilon}_n\triangleq
z^2\underline{s}^2(z)\psiup_n(z)+\bar{\epsilon}_n .
\]

By the identity of quadric form \eqref{quadric} and the fact, from
Lemma \ref{diagonal}, that $\mathbb{E}[D_j^{-1}(z)]_{ii}$ can be
replaced by
$s(z)=-z^{-1}(\underline{s}(z)+1)^{-1}$, we have
\begin{eqnarray} \label{sum}
\bar{J}_{n1}(z)
&=&-z^2\underline{s}^2(z)\sum_{j=1}^n\mathbb{E}\varepsilon^2_j(z)+\bar{\epsilon}_n\nonumber\\
&=&-\frac{z^2\underline{s}^2(z)}{n^2}\sum_{j=1}^n\mathbb{E}
\Biggl[\sum_{i=1}^p\kappa_2[D_j^{-1}(z)]^2_{ii}+\kappa
_1\operatorname{tr}D_j^{-2}(z)+\operatorname{tr}D_j^{-2}(z) \Biggr]+\bar{\epsilon}_n\\
\nonumber&=&y_n\kappa_2k^2(z)-z^2\underline{s}^2(z)(\kappa_1+1)\psiup
_n(z)+\bar{\epsilon}_n ,
\end{eqnarray}
where $\kappa_1$, $\kappa_2$ and $k(z)$ were defined in Theorem \ref{main}. Our goal is now to find the limit of $\psiup_n(z)$. Using the
expansion of $D_j^{-1}(z)$ in \eqref{djni}, we get
\begin{eqnarray*}
\psiup_n(z)&=&\frac{1}{n^2}\sum_{j=1}^n\frac{p}{(z+z\underline{s}(z))^2}
+z^2\underline{s}^2(z)\frac{1}{n^2}\sum_{j=1}^n\mathbb{E}n\operatorname{tr}A^2(z)+\bar{\epsilon}_n\\
&=&\frac{k^2(z)}{n^2}\sum_{j=1}^n\sum^n_{i,l\neq j}
\mathbb{E}\operatorname{tr} \biggl[ \biggl(r_ir_i^*-\frac{1}{n}I
\biggr)D_{ij}^{-1}(z)D_{lj}^{-1}(z) \biggl(r_lr_l^*-\frac{1}{n}I \biggr) \biggr]\\
&&{}+\frac{1}{n^2}\sum_{j=1}^n\frac{p}{z^2(\underline{s}(z)+1)^2}+\bar{\epsilon}_n .
\end{eqnarray*}
Note that the cross terms will be $0$ if either $D_{ij}^{-1}(z)$ or
$D_{lj}^{-1}(z)$ is replaced by $D_{lij}^{-1}(z)$, where
$D_{lij}(z)=D_{ij}(z)-r_lr_l^*=D_{lj}^{-1}(z)-r_ir_i^*$ and
\[
D_{ij}^{-1}(z)-D_{lij}^{-1}(z)=-\frac
{D_{lij}^{-1}(z)r_lr_l^*D_{lij}^{-1}(z)}{1+r_l^*D_{lij}^{-1}(z)r_l}.
\]

Therefore, by \eqref{d-dj}, we conclude that the sum of cross terms is
negligible and bounded by $K/(\sqrt{n}v^3)$. Thus, we find that
\begin{eqnarray*}
&&\frac{1}{n^2}\sum_{j=1}^n\sum^n_{i,l\neq j}
\mathbb{E}\operatorname{tr} \biggl[ \biggl(r_ir_i^*-\frac{1}{n}I
\biggr)D_{ij}^{-1}(z)D_{lj}^{-1}(z) \biggl(r_lr_l^*-\frac{1}{n}I \biggr) \biggr]\\
&&\quad=\frac{1}{n^2}\sum_{j=1}^n\sum^n_{i\neq j}
\mathbb{E}\operatorname{tr} \biggl[ \biggl(r_ir_i^*-\frac{1}{n}I \biggr)D_{ij}^{-2}(z)
\biggl(r_ir_i^*-\frac{1}{n}I \biggr) \biggr]+\bar{\epsilon}_n\\
&&\quad=\frac{1}{n^2}\sum_{j=1}^n\sum^n_{i\neq j}
\mathbb{E} [(r_i^*D_{ij}^{-2}(z)r_i)(r_i^*r_i) ]+\bar{\epsilon}_n\\
&&\quad=\frac{1}{n^2}\sum_{j=1}^n\sum^n_{i\neq j}
\frac{1}{n^2}\mathbb{E} \bigl[\operatorname{tr}D_{ij}^{-2}(z)\bigl(p+\mathrm{O}(1)\bigr)
\bigr]+\bar{\epsilon}_n=y_n\psiup_n(z)+\bar{\epsilon}_n .
\end{eqnarray*}

From above, we get that
\[
\psiup_n(z)=\frac{y_n}{z^2(\underline{s}(z)+1)^2}
+y_nk^2(z)\psiup_n(z)+\bar{\epsilon}_n .
\]

Combined with \eqref{sum}, we have
\[
J_n(z)=\kappa_2y_nk^2(z)-\frac{\kappa_1y_nk^2(z)}{1-y_nk^2(z)}+\bar{\epsilon}_n .
\]

Thus, from \eqref{diff}, it follows that
\begin{eqnarray*}
M_n(z)&=&n\mathbb{E}\underline{s}_n(z)\underline{s}_n^0(z)R_n(z)\big/\biggl(1-\frac
{y_n\mathbb{E}\underline{s}_n(z)\underline{s}_n^0(z)}{(\underline
{s}_n^0(z)+1)(\mathbb{E}\underline{s}_n(z)+1)}\biggr)\\
&=&-\frac{\underline{s}_n^0(z)}{\mathbb{E}\underline{s}_n(z)+1}J_n(z)\big/
\biggl(1-\frac{y_n\mathbb{E}\underline{s}_n(z)\underline
{s}_n^0(z)}{(\underline{s}_n^0(z)+1)(\mathbb{E}\underline{s}_n(z)+1)}\biggr)\\
&=&\frac{\kappa_1y_nk^3(z)}{(1-y_nk^2(z))^2}
-\frac{\kappa_2y_nk^3(z)}{1-y_nk^2(z)}
+\bar{\epsilon}_n\\
&\triangleq&\widetilde{M}_1(z)+ \widetilde{M}_2(z)+\bar{\epsilon}_n .
\end{eqnarray*}

Therefore, we can calculate the mean function in the following two parts:
\begin{eqnarray*}
&&-\frac{1}{2\curpi \mathrm{i}}\int_{\gamma_{mh}}f_m(z)\widetilde{M}_1(z)\,\mathrm{d}z\\
&&\quad =-\frac
{\kappa_1}{2\curpi \mathrm{i}}\int_{\gamma_{mh}}f_m(z)\frac
{y_nk^3(z)}{(1-y_nk^2(z))^2}\,\mathrm{d}z\\
&&\quad=\frac{\kappa_1}{4\curpi \mathrm{i}}\int_{\gamma_{mh}}f_m(z)\frac{\mathrm{d}}{\mathrm{d}z}\ln
\bigl(1-y_nk^2(z)\bigr)\,\mathrm{d}z =-\frac{\kappa_1}{4\curpi \mathrm{i}}\int_{\gamma_{mh}}f'_m(z)\ln
\bigl(1-y_nk^2(z)\bigr)\,\mathrm{d}z\\
&&\quad\longrightarrow\frac{\kappa_1}{2\curpi}\int^b_a f'(x)\operatorname{arg}
\bigl(1-yk^2(x) \bigr)\,\mathrm{d}x ,
\end{eqnarray*}
as $n\rightarrow\infty$, $a_{\mathrm{l}}\rightarrow a$ and $b_{\mathrm{r}}\rightarrow b$; similarly,
\begin{eqnarray*}
-\frac{1}{2\curpi \mathrm{i}}\int_{\gamma_{mh}}f_m(z)\widetilde{M}_2(z)\,\mathrm{d}z
&=& \frac{\kappa_2}{2\curpi \mathrm{i}}\int_{\gamma_{mh}}f_m(z)\frac
{y_nk^3(z)}{1-y_nk^2(z)}\,\mathrm{d}z\\
&\longrightarrow&-\frac{\kappa_2}{\curpi}\int^b_af(x)\Im\biggl(\frac
{yk^3(x)}{1-yk^2(x)} \biggr)\,\mathrm{d}x .
\end{eqnarray*}

Hence, summing the two terms, we obtain the mean function of the
limiting distribution in~\eqref{mean}.

\begin{appendix}
\section*{Appendix}\label{app}
\setcounter{section}{1}
\setcounter{equation}{0}

\begin{lem} \label{cr}
Under the conditions in Theorem \ref{main}, we have
\begin{eqnarray*}
\|\mathbb{E}F_n-F\|&=&\mathrm{O}(n^{-1/2}) ,\qquad \|F_n-F\|=\mathrm{O}_p(n^{-2/5}) , \\
\| F_n-F\|&=&\mathrm{O}(n^{-2/5+\eta}) \qquad\mbox{a.s. for any }  \eta>0 .
\end{eqnarray*}
\end{lem}

This follows from Theorems 1.1, 1.2 and 1.3 in \cite{bmy2003}.

\begin{lem}[{[Burkholder (1973), \cite{burk1973}]}] \label{Burk}
Let $X_k$, $k=1,2,\ldots,$ be a complex martingale difference sequence
with respect to the increasing $\sigma$-fields $\mathcal{F}_k$. Then,
for $p>1$,
\[
\mathbb{E} \Bigl|\sum X_k \Bigr|^p\leq K_p\mathbb{E} \Bigl(\sum|X_k|^2 \Bigr)^{p/2}.
\]
\end{lem}

In the reference \cite{burk1973}, only real variables were considered.
It is straightforward to extend to complex cases.

%
\begin{lem} \label{quadric1}
For $x=(x_1,\ldots,x_n)^{\mathrm{t}}$ with i.i.d. standardized real or complex
entries such that
$\mathbb{E}x_i=0$ and $\mathbb{E}|x_i|^2=1$, and for C an $n\times
n$ complex matrix, we have, for any $p\geq2,$
\[
\mathbb{E}|x^*Cx-\operatorname{tr}C|^p\leq K_p [(\mathbb
{E}|x_1|^4\operatorname{tr}CC^*)^{p/2}+\mathbb{E}|x_1|^{2p}\operatorname{tr}(CC^*)^{p/2} ].
\]
\end{lem}

This is Lemma 8.10 in \cite{bai06}.

%
\begin{lem}\label{quadric2}
For any non-random $p\times p$ matrix $A$,
\[
\mathbb{E}|r_1^*Ar_1|^2\leq Kn^{-1}\|A\|^2 .
\]
\end{lem}

\begin{pf}
For non-random $p\times p$ matrix $A$,
\begin{eqnarray*}
\mathbb{E}|r_1^*Ar_1|^2 &=&\frac{1}{n^2}\mathbb{E} \Biggl|\sum_{l,k=1}^p\bar
{x}_{l1}a_{lk}x_{k1} \Biggr|^2\\
&=&\frac{1}{n^2}\mathbb{E} \Biggl(\sum_{l\neq
k}^p\bar{x}_{l1}^2a_{lk}^2x_{k1}^2+\sum_{l\neq
k}^p|{x_{l1}|^2|x_{k1}|^2a_{lk}a_{kl}}+\sum_{l=1}^p|x_{1l}|^4a_{ll}^2
\Biggr)\\
&\leq&\frac{K}{n^2}\mathbb{E} \Biggl(\sum_{l,k=1}^p|a_{lk}|^2
\Biggr)=Kn^{-2}\mathbb{E}\operatorname{tr}(A\bar{A})\leq Kn^{-1}\|A\|^2 .
\end{eqnarray*}
\upqed
\end{pf}

\begin{lem}\label{quadric3}
For non-random $p\times p$ matrices $A_k,k=1,\ldots,s,$
\begin{equation}\label{6.1}
\mathbb{E} \Biggl|\prod_{k=1}^s \biggl(r_1^*A_lr_1-\frac{1}{n}\operatorname{tr}A_l \biggr) \Biggr|\leq
Kn^{-((s/2)\wedge3)}\delta_n^{2(s-4)\vee0}\prod_{k=1}^s\|A_l\| .
\end{equation}
\end{lem}

\begin{pf}
Recalling the truncation steps $\mathbb{E}|x_{11}|^8<\infty$ and Lemma
\ref{quadric1}, we have, for all $l>1$,
\begin{eqnarray}\label{6.2}
\mathbb{E}|r_1^*A_1r_1-n^{-1}\operatorname{tr}A_1|^l&\leq& K\|A_1\|
^ln^{-l}\bigl(n^{l/2}+\bigl(\sqrt{n}\delta_n\bigr)^{(2l-8)\vee0}n\bigr)\nonumber\\[-8pt]\\[-8pt]
&=&K\|A_1\|^ln^{-((l/2)\wedge3)}(\delta_n)^{2(l-4)\vee0} .\nonumber
\end{eqnarray}
Then, \eqref{6.1} is the consequence of \eqref{6.2} and the H\"{o}lder
inequality.
\end{pf}

\begin{lem}\label{beta}
Under the conditions in Theorem \ref{main}, for any $l\geq2,   \mathbb{E}|\beta
_j(z)|^l$, $\mathbb{E}|\bar{\beta}_j(z)|^l$ and $|b_n(z)|^l$ are uniformly
bounded in $\gamma_{mh}$. Furthermore, $\beta_j(z),  \bar{\beta}_j(z)$
and $b_n(z)$ are uniformly convergent in probability to $-z\underline
{s}(z)$ in
$\gamma_{mh}$.
\end{lem}

\begin{pf}
By (4.2) and (4.3) in \cite{bs1998}, we have, for any $l\leq2$,
\begin{eqnarray} \label{63}
\mathbb{E}|\operatorname{tr}D_j^{-1}(z)-\mathbb{E}\operatorname{tr}D_j^{-1}(z)|^l&\leq& Kn^{
{l}/{2}}v^{-l} ,\\
\label{64}
\mathbb{E}|r_jD_j^{-1}(z)r_j-1/n\mathbb
{E}\operatorname{tr}D_j^{-1}(z)|^l&\leq& Kn^{-{l}/{2}}v^{-l} .
\end{eqnarray}
This lemma follows from Lemma \ref{quadric1}, \eqref{63}, \eqref{64}
and the following facts.

\begin{fact}
Since $s_n^0(z)=-\frac{1}{2} (\frac{1}{y_n}-\frac{1}{y_nz}\sqrt
{z^2-(1+y_n)z+(1-y_n)^2}-\frac{1-y_n}{y_nz} )$ and
$\underline{s}_n^0(z)=-\frac{1-y_n}{z}+y_ns_n^0(z)$, we have
\[
z\underline{s}_n^0(z)=-\tfrac{1}{2} \bigl(1-y_n+z-\sqrt
{z^2-(1+y_n)z+(1-y_n)^2} \bigr) .
\]

Thus, $z\underline{s}_n^0(z)$ is bounded in any bounded and closed
complex region.
\end{fact}

\begin{fact}
\begin{eqnarray*}
|b_n(z)-\mathbb{E}\beta_j(z)| &\leq&\frac{1}{v^2}\mathbb{E}
\biggl|r_j^*D_j^{-1}(z)r_j-\frac{1}{n}\mathbb{E}\operatorname{tr}D_j^{-1}(z) \biggr|\\
&\leq&\frac{1}{v^2} \biggl[\mathbb{E}|\varepsilon_j(z)|+
\mathbb{E} \biggl|\frac{1}{n}\operatorname{tr}D_j^{-1}(z)-\frac{1}{n}\mathbb
{E}\operatorname{tr}D_j^{-1}(z) \biggr| \biggr]\\
&\leq&\frac{1}{v^2} \biggl[\frac{K}{\sqrt{n}v}+\frac{K}{\sqrt{n}v} \biggr]=\frac
{K}{\sqrt{n}v^3} ,
\end{eqnarray*}
where the last inequality follows from \eqref{ded}.
\end{fact}

\begin{fact}
Taking expectation on \eqref{snbeta}, one can find
\[
z\mathbb{E}\underline{s}_n(z)=-\frac{1}{n}\sum_{j=1}^n\mathbb{E}\beta
_j(z)=-\mathbb{E}\beta_j(z) .
\]
\end{fact}

\begin{fact}
From Lemma \textup{\ref{cr}}, we have
\begin{eqnarray*}
|z\mathbb{E}\underline{s}_n(z)-z\underline{s}_n^0(z)| &\leq&
zy_n\mathbb{E}|s_n(z)-s_n^0(z)|\\
&=&zy_n\mathbb{E} \biggl|\int\frac
{1}{x-z}(F^{B_n}-F^{y_n})(\mathrm{d}x) \biggr|\\
&\leq&
\frac{K}{v}\|F^{B_n}-F^{y_n}\|\\
&=&\frac{K}{v}\mathrm{O}_p(n^{-
{2}/{5}})=\mathrm{O}_p(n^{-{2}/{5}}v^{-1}) .
\end{eqnarray*}
\end{fact}
\upqed
\end{pf}

\begin{lem}\label{diagonal}
Under the conditions in Theorem \ref{main}, as $n\rightarrow\infty$,
\[
\max_{i,j}|\mathbb{E}_j[D_j^{-1}(z)]_{ii}-s(z)|\rightarrow
0  \qquad \mbox{in probability}
\]
uniformly in $\gamma_{mh}$, where the maximum is taken over all $1\leq
i\leq p$ and $1\leq j\leq n$.
\end{lem}

\begin{pf}
First, let $e_j$ $(1\leq j\leq n)$ be the $p$-vector whose $j$th
element is 1, the rest being 0 and $e'_i$, the transpose of $e_i$.
Then,
\begin{eqnarray*}
\mathbb{E}|[D^{-1}(z)]_{ii}-[D_j^{-1}(z)]_{ii}|
&=&\mathbb{E}\bigl|e'_i\bigl(D^{-1}(z)-D^{-1}_j(z)\bigr)e_i\bigr|\\
&=&\mathbb{E}|\beta_j(z)e'_iD^{-1}_j(z)r_jr^*_jD^{-1}_j(z)e_i|\\
&\leq&(\mathbb{E}|\beta_j(z)|^2)^{1/2}(\mathbb
{E}|r^*_jD^{-1}_j(z)e_ie'_iD^{-1}_j(z)r_j|^2)^{1/2}
\leq\frac{K}{\sqrt{n}v^2} .
\end{eqnarray*}

Second, by martingale inequality, for any $\epsilon>0$, we have
\begin{eqnarray*}
&&P\Bigl(\max_{i,j}|\mathbb{E}_j[D^{-1}(z)]_{ii}-\mathbb
{E}[D^{-1}(z)]_{ii}|>\epsilon\Bigr)\\
&&\quad\leq\sum_{i=1}^pP\Bigl(\max_{j}|\mathbb
{E}_j[D^{-1}(z)]_{ii}-\mathbb{E}[D^{-1}(z)]_{ii}|>\epsilon\Bigr)\\
&&\quad\leq \sum_{i=1}^p\frac{1}{\epsilon^6}\mathbb
{E}|[D^{-1}(z)]_{ii}-\mathbb{E}[D^{-1}(z)]_{ii}|^6\\
&&\quad= \frac{1}{\epsilon^6}\sum_{i=1}^p\mathbb{E} \Biggl|\sum_{l=1}^n(\mathbb
{E}_l-\mathbb{E}_{l-1})
\beta_l(z)e'_iD_l^{-1}(z)r_lr^*_lD_l^{-1}(z)e_i \Biggr|^6\\
&&\quad\leq K\sum_{i=1}^p\mathbb{E} \Biggl(\sum_{l=1}^n|(\mathbb{E}_l-\mathbb{E}_{l-1})
\beta_l(z)e'_iD_l^{-1}(z)r_lr^*_lD_l^{-1}(z)e_i|^2 \Biggr)^3 .
\end{eqnarray*}
Let $Z_l(z)=e'_iD_l^{-1}(z)r_lr^*_lD_l^{-1}(z)e_i$. We have that
\[
|\mathbb{E}Z_l(z)|\leq\frac{K}{nv^2}\quad\mbox{and}\quad
\mathbb{E}|Z_l(z)-\mathbb{E}Z_l(z)|^2\leq\frac{K}{n^2v^4} .
\]
Thus, we obtain
\begin{eqnarray*}
&&P\Bigl(\max_{i,j}|\mathbb{E}_j[D^{-1}(z)]_{ii}-\mathbb
{E}[D^{-1}(z)]_{ii}|>\epsilon\Bigr)\\
&&\quad\leq K\sum_{i=1}^p\mathbb{E} \Biggl(\sum_{l=1}^n\frac{K}{n^2v^4} \Biggr)^3=\frac
{K}{n^2v^{12}} .
\end{eqnarray*}

Finally,
\[
\mathbb{E}[D^{-1}]_{ii}=\frac{1}{p}\sum_{i=1}^p\mathbb{E}[D^{-1}]_{ii}=\mathbb{E}s_n(z) .
\]

In Section \ref{mf}, it is proved that $p(\mathbb{E}s_n(z)-s(z))$ converges
to 0 uniformly on $\gamma_{mh}$. The proof of Lemma \ref{diagonal} is thus complete.
\end{pf}
\end{appendix}

\section*{Acknowledgements}
The authors would like to thank one Associate Editor and one
referee for many constructive comments. The first two authors were
partially supported by NSFC 10871036 and NUS Grant R-155-000-079-112.
Wang Zhou was partially supported by Grant
R-155-000-083-112 at
the National University of Singapore.

\printhistory

\end{document}